\documentclass[reqno,10pt]{amsart}
\usepackage{amsmath,amsthm,amssymb,amsfonts,color,xcolor,pifont,cite}
\usepackage{a4wide}
\usepackage[colorlinks]{hyperref}
\usepackage[showonlyrefs]{mathtools}
\mathtoolsset{showonlyrefs=true}
\usepackage{framed,enumerate}
\usepackage[normalem]{ulem}
\usepackage{mathrsfs}

\newtheorem{theorem}{Theorem}[section]
\newtheorem{proposition}[theorem]{Proposition}
\newtheorem{corollary}[theorem]{Corollary}
\newtheorem{lemma}[theorem]{Lemma}

\newtheorem{remark}[theorem]{Remark}

\newcommand{\R}{\mathbb R}

\def\eps{\varepsilon}
\def\les{\lesssim}

\def\ges{\gtrsim}
\author{Michael Goldman}
\address{LAMA, CNRS, Universit\'e Paris-Est Cr\'eteil, 61 avenue du G\'en\'eral de Gaulle, 94010 Cr\'eteil,
	France}
\email{michael.goldman@cnrs.fr}

\author{Matteo Novaga}
\address{Department of Mathematics, University of Pisa, 56127 Pisa, Italy}
\email{matteo.novaga@unipi.it}

\author{Berardo Ruffini}
\address{Department of Mathematics, University of Bologna, 40126 Bologna, Italy}
\email{berardo.ruffini@unibo.it}

\numberwithin{equation}{section}

\title[Sharp quantitative fractional Faber--Krahn]{The fractional
Faber--Krahn inequality in sharp quantitative form}


\begin{document}

\begin{abstract}
We prove that the first Dirichlet
eigenvalue of the fractional Laplacian satisfies the
sharp quantitative Faber--Krahn inequality. This was conjectured by Brasco,
Cinti and Vita, who proved the inequality with a sub-optimal exponent.  The
proof adapts to the non-local setting the level-set argument recently
introduced by Guerra, Machado and Ramos for the Laplacian. It is carried
out on the horizontal level sets of the Caffarelli--Silvestre extension
of the torsion function, where the asymmetry  appears only as
a boundary value at the corner of the extension variables.
\end{abstract}

\maketitle

\section{Introduction}

Let $n\geq2$, $s\in(0,1)$, and let $\Omega\subset\R^n$ be an open set
with $0<|\Omega|<\infty$. We denote by $\lambda_s(\Omega)$ the first
eigenvalue of the fractional Laplacian $(-\Delta)^s$ in $\Omega$ with
homogeneous exterior Dirichlet condition; the precise normalization is
fixed in Section~\ref{sec:preliminaries}. If $B_\Omega$ is a ball with
$|B_\Omega|=|\Omega|$, the fractional Faber--Krahn inequality states
that $\lambda_s(\Omega)\geq\lambda_s(B_\Omega)$; it follows from the
fractional P\'olya--Szeg\H{o} inequality, see \cite{AL, BCV}. A relatively natural question in this context is
 to quantify this inequality and prove that sets with $\lambda_s(\Omega)$ close to $\lambda_s(B_\Omega)$ are close to $B_\Omega$, up to translation.
  See e.g. \cite{BDV,BCV,GMR,FJ,FMP,BDP}.
  We measure
the distance of $\Omega$ from balls by its Fraenkel asymmetry
\[
 \mathcal A(\Omega)=\inf_{a\in\R^n}
 \frac{|\Omega\mathbin\triangle B_\Omega(a)|}{|\Omega|},
\]
where $B_\Omega(a)$ is the ball of center $a$ with $|B_\Omega(a)|=|\Omega|$.

Brasco, Cinti and Vita \cite[Theorem~1.3]{BCV}, proved that for every open set of finite measure

\[\lambda_s(\Omega)/\lambda_s(B_\Omega)-1\geq c\,\mathcal A(\Omega)^{3/s}\]
 and conjectured
\cite[Remark~1.6]{BCV} that the exponent $3/s$ can be replaced by $2$. For the torsional rigidity, the exponent was improved to $2+1/s$ in \cite[Theorem~1.8]{Kim}.
The exponent $2$ is optimal also in the fractional setting, as shown by the volume-preserving ellipsoidal perturbations in Lemma~\ref{lem:ellipsoidal-sharpness}.
 Our main
result confirms this conjecture.

\begin{theorem}\label{thm:main-fk}
For every $n\geq2$ and $s\in(0,1)$ there is $c=c(n,s)>0$ such that
every open set $\Omega\subset\R^n$ with $0<|\Omega|<\infty$ satisfies
\begin{equation}\label{eq:target-fk}
 \frac{\lambda_s(\Omega)}{\lambda_s(B_\Omega)}-1
 \geq c\,\mathcal A(\Omega)^2.
\end{equation}
\end{theorem}
As in \cite{BDV,GMR}, we actually prove the corresponding quantitative 
inequality for the fractional torsional rigidity.
 Let
$T_s(\Omega)$ denote the fractional torsional rigidity of $\Omega$ (see
Section~\ref{sec:preliminaries}), and let $B=B_1$ be the unit ball.

\begin{theorem}\label{thm:main-sv}
For every $n\geq2$ and $s\in(0,1)$ there is $C=C(n,s)>0$ such that
every open set $\Omega\subset\R^n$ with $|\Omega|=|B|$ satisfies
\begin{equation}\label{eq:main-sv}
 \mathcal A(\Omega)^2\leq C\bigl(T_s(B)-T_s(\Omega)\bigr).
\end{equation}
\end{theorem}

Theorem~\ref{thm:main-sv} is proved in Section~\ref{sec:reduction}. Theorem~\ref{thm:main-fk}
follows from it and from the fractional K\"ohler--Jobin inequality from \cite{BDFPV}.
 Since the argument is identical to that in the local case \cite{BDV,GMR}, we omit the proof of this last derivation.

As in \cite{GMR}, the constants in Theorems~\ref{thm:main-fk} and
\ref{thm:main-sv} are in principle computable. We do not track their
numerical values. We expect the proof to apply mutatis mutandis to the case of the fractional capacity, which would lead to
 a similar sharp extension of the results obtained in \cite{cinti2021quantitative}.

\subsection*{Strategy of the proof}

Our proof of Theorem~\ref{thm:main-sv} follows the strategy of \cite{GMR}. 
We refine the level-set/symmetrization  approach of \cite{BCV} in a similar spirit as \cite{GMR} 
  refines  the proof from \cite{fusco2009stability}. For $\Omega$ a given set, we write $u_\Omega$ for its torsion function and $u^*_\Omega$ for the corresponding Schwarz rearrangement.
  Although not always explicitly presented in this form, the starting point in \cite{fusco2009stability,BCV} is to split the torsion deficit as the sum of a deficit between the energy of $u_\Omega$ and $u_\Omega^*$ plus
  the energy of $u_\Omega^*-u_B$, see Lemma \ref{lem:split} here. By coarea, this is then transferred to a control on the isoperimetric deficit 
  of the  level sets of $u_\Omega$. It is then natural to try to use the quantitative isoperimetric inequality from \cite{FMP} to control the Fraenkel asymmetry of these level sets. 
  The difficulty is then to relate these asymmetries to the asymmetry of $\Omega$. In \cite{fusco2009stability,BCV}, this is done by a mean-value-type argument which leads to a nonsharp exponent.
   The central idea of \cite{GMR} is to replace this Chebyshev-type argument by the fundamental theorem of calculus which then requires estimating the derivative
    with respect to the height, or somewhat equivalently with respect to volume, of this Fraenkel asymmetry. 
    Interestingly, bounding this derivative involves the strong form of the quantitative isoperimetric inequality, due to Fusco and Julin \cite{FJ} (see also \cite{hensel2026quantitative}) 
     rather than the one involving the Fraenkel asymmetry from \cite{FMP}.

As in \cite{BCV}, we work with the Caffarelli--Silvestre \cite{CS} extension $U$ of
the torsion function $u_\Omega$  and with the horizontal level
sets $E_{t,z}=\{U(\cdot,z)>t\}$ at each height $z>0$. Rearranging $U$ in
the horizontal variables, we write the torsion deficit as the weighted
integral of four nonnegative defects (Lemma~\ref{lem:half-four-deficits}).
The first two defects $D_H$ and $D_I$ are the same as in \cite{GMR}:
$D_H$ comes from the Cauchy--Schwarz inequality on
$\partial E_{t,z}$ and measures the oscillation of $|\nabla_xU|$, while
$D_I$ is the isoperimetric defect of $E_{t,z}$.
The last two defects are new compared with the local case. The
vertical defect $D_V$ measures the oscillation of $\partial_zU$ on
$\partial E_{t,z}$, and $D_B$ measures the weighted distance between the
horizontally rearranged extension and the extension of the torsion
function of the ball.

We normalize $|\Omega|=\omega_n$. Following \cite{GMR}, at every height
$z$ we parametrize the level sets by the radius of a ball of the same
volume rather than by $t$. Thus we choose $t(\rho,z)$ so that
$|E_{t(\rho,z),z}|=\omega_n\rho^n$, set
$E_{\rho,z}=E_{t(\rho,z),z}$, and consider the unnormalized Fraenkel
asymmetry, or moving ball distance in the language of \cite{GMR},
$q(\rho,z)=\inf_a|E_{\rho,z}\triangle B_\rho(a)|$.
  To follow the proof of \cite{GMR} closely, we would ideally like to control the derivative of $q$ along a curve connecting the point $(1,0)$ to a point $(\rho,z)$ with $\rho,z\simeq 1/2$. The corner value
  $q(1,0)=\omega_n\mathcal A(\Omega)$, understood as a joint limit, is thus the quantity we want to control. However, there are two difficulties. The first, which might be only technical, concerns regularity in $z$. Following \cite{GMR}, we can prove that $q$ is Lipschitz, hence absolutely continuous, in the horizontal variable $\rho$, see \eqref{eq:qlip}, but we are not able to obtain a similar result in the $z$ variable.
    The key insight, see Proposition \ref{prop:half-vertical-flux} and Remark \ref{rem:vertic-deriv}, is that while $q$ is not necessarily absolutely continuous in $z$,   the weighted quantity 
     $wq/w_B$ is better behaved. Here $w$ is the radial speed of the level-sets of the extension of $u_\Omega$ (and similarly for $w_B$). The second difficulty is that the weights are degenerating
      as $(\rho,z)$ approaches $(1,0)$. To circumvent this issue, we replace pointwise bounds by averaged ones on horizontal windows of variable size in the proof of Theorem~\ref{thm:half-confined-SV}.

As in \cite{BDV,GMR}, we actually first prove Theorem \ref{thm:main-sv} for sets which have a diameter bounded by a universal constant. 
To reduce the general case to this one we use an iterated truncation/De Giorgi type argument in Proposition \ref{prop:half-truncation}   in a similar way as in \cite[Lemma~5.3]{BDV} and
\cite[Lemma~4.1]{GMR}. 

\subsection*{Organization of the paper}
Section~\ref{sec:preliminaries} fixes the notation and proves the
splitting of the torsion deficit into a horizontal rearrangement loss and
a distance to the ball. Section~\ref{sec:bounded-torsion} proves the
sharp quantitative Saint-Venant inequality for sets contained in a fixed
ball. Section~\ref{sec:reduction} removes the confinement and
Section~\ref{sec:eigenvalue} proves the sharpness of the exponent.
Appendix~\ref{app:ball-hardy} contains the estimates for the extension of
the torsion function of the ball and the boundary control of the radial
profile.
 \subsection*{AI disclosure and development of this project} We started thinking about the question of
  the quantitative stability for the fractional capacity after completing the work \cite{goldman2024rigidity}. After
   realizing that the approach there would most likely not yield the desired result, we started exploring the adaptation to
    the non-local setting of the selection principle strategy from \cite{BDV}. We faced, however, several technical difficulties,
     the first of which is the lack of computations of second-order shape derivatives for non-local problems, see \cite{DFW,DGV,muratov2022conducting} for (tedious) 
     first-order computations. In spring 2026, with the help of LLMs, we made substantial progress on that 
     research program and we are reasonably optimistic that
      it could be carried out at the price of extensive technical work. 
      As an illustration of what this route might look like, we provide in the companion paper \cite{gnrsecond} the computation of
       the second-order shape derivative of the torsion in the case of smooth sets, and show linear quantitative stability of the ball. 
       After the paper \cite{GMR} appeared on arXiv it was very tempting to check if their much more direct strategy could
        be adapted to the  non-local setting. This was done principally by the LLMs with some interactions and simplifications from our end.  
        The first draft of the paper was also generated before thorough human inspection and editing. We take
         full responsibility for the content of this paper.  

\section{Preliminaries and notation}\label{sec:preliminaries}

Fix $n\geq2$ and $s\in(0,1)$. Constants denoted by $C,c>0$ may
depend on $n,s$; dependence on any further parameter is displayed. We write  $A \les B$  for $A \le C B$ and $A\asymp B$ if $ A\les  B \les  A$. For
an open set $\Omega\subset\R^n$ of finite measure, set
\[
 H_0^s(\Omega)=\overline{C_c^\infty(\Omega)}^{\,H^s(\R^n)}
\]
and
\[
 \mathcal E_s(u,v)=\frac{c_{n,s}}2
 \iint_{\R^n\times\R^n}
 \frac{(u(x)-u(y))(v(x)-v(y))}{|x-y|^{n+2s}}\,dx\,dy,
 \qquad
 c_{n,s}=\frac{4^s\Gamma(n/2+s)}{\pi^{n/2}|\Gamma(-s)|}.
\]
We write $\mathcal E_s(u)=\mathcal E_s(u,u)$. With the Fourier
convention
$\widehat f(\xi)=\int_{\R^n}e^{-ix\cdot\xi}f(x)\,dx$, this is the
quadratic form associated with the multiplier $|\xi|^{2s}$. Thus
\[
 (-\Delta)^sf(x)=c_{n,s}\operatorname{PV}\!\int_{\R^n}
 \frac{f(x)-f(y)}{|x-y|^{n+2s}}\,dy
\]
whenever the principal value is defined. The first eigenvalue is
\[
 \lambda_s(\Omega)=\inf_{0\ne v\in H_0^s(\Omega)}
 \frac{\mathcal E_s(v)}{\|v\|_2^2}.
\]

The torsion function $u_\Omega$ is the unique minimizer of
\begin{equation}\label{eq:torsion-variational}
 u_\Omega=\operatorname*{argmin}_{v\in H_0^s(\Omega)}
 \left\{\frac12\mathcal E_s(v)-\int_\Omega v\,dx\right\}.
\end{equation}
Equivalently, for every $\varphi\in H_0^s(\Omega)$,
\begin{equation}\label{eq:torsion-weak}
 \mathcal E_s(u_\Omega,\varphi)=\int_\Omega\varphi\,dx,
\end{equation}
which is the weak formulation of
\[
 \begin{cases}
 (-\Delta)^su_\Omega=1&\text{in }\Omega,\\
 u_\Omega=0&\text{in }\R^n\setminus\Omega.
 \end{cases}
\]
The torsional rigidity is
\begin{equation}\label{eq:torsion-rigidity}
 T_s(\Omega)=\int_\Omega u_\Omega\,dx=\mathcal E_s(u_\Omega)
 =\max_{v\in H_0^s(\Omega)}
 \left\{2\int_\Omega v\,dx-\mathcal E_s(v)\right\}.
\end{equation}

Since $u_{r\Omega}(x)=r^{2s}u_\Omega(x/r)$ for $r>0$, we have
\begin{equation}\label{eq:scaling}
 T_s(r\Omega)=r^{n+2s}T_s(\Omega),\qquad
 \lambda_s(r\Omega)=r^{-2s}\lambda_s(\Omega).
\end{equation}
We work under the normalization
$|\Omega|=\omega_n=|B|$, where $B=B_1$, and set
\begin{equation}\label{eq:absolute-torsion-deficit}
 \delta_T(\Omega)=T_s(B)-T_s(\Omega).
\end{equation}

The Caffarelli--Silvestre extension of $f\in H^s(\R^n)$ is
\begin{equation}\label{eq:general-poisson}
 \operatorname{Ext}_sf(x,z)=P_z^{(s)}*f(x),\qquad
 P_z^{(s)}(x)=\frac{\Gamma(n/2+s)}{\pi^{n/2}\Gamma(s)}
 \frac{z^{2s}}{(|x|^2+z^2)^{(n+2s)/2}}.
\end{equation}
It solves $\operatorname{div}(z^{1-2s}\nabla\operatorname{Ext}_sf)=0$
in $\R^{n+1}_+=\R^n\times(0,\infty)$ and \cite{CS}
\begin{equation}\label{eq:general-extension-normalization}
 d_s=2^{2s-1}\frac{\Gamma(s)}{\Gamma(1-s)},\qquad
 \mathcal E_s(f)=d_s\int_{\R^{n+1}_+}z^{1-2s}
 |\nabla\operatorname{Ext}_sf|^2\,dx\,dz.
\end{equation}
Moreover,
\begin{equation}\label{varextension}
  \mathcal E_s(f)=d_s\min_{F(x,0)=f}\int_{\R^{n+1}_+}z^{1-2s}
 |\nabla F|^2\,dx\,dz
\end{equation}
and 
\[
 (-\Delta)^sf=-d_s\lim_{z\downarrow0}z^{1-2s}
 \partial_z\operatorname{Ext}_sf
\]
in the distributional trace sense. Hence $U=\operatorname{Ext}_su_\Omega$
satisfies
\[
 \begin{cases}
 \operatorname{div}(z^{1-2s}\nabla U)=0&\text{in }\R^{n+1}_+,\\
 U=0&\text{on }(\R^n\setminus\Omega)\times\{0\},\\
 -d_s\displaystyle\lim_{z\downarrow0}z^{1-2s}U_z=1
 &\text{on }\Omega\times\{0\}.
 \end{cases}
\]
By this we mean that  for every finite-energy $\Phi$ whose trace belongs to
$H_0^s(\Omega)$,
\begin{equation}\label{eq:torsion-extension-weak}
 d_s\int_{\R^{n+1}_+}z^{1-2s}\nabla U\cdot\nabla\Phi\,dx\,dz
 =\int_\Omega\Phi(x,0)\,dx.
\end{equation}

For a nonnegative measurable function $f$ vanishing at infinity we
denote by $f^*$ its Schwarz symmetric decreasing rearrangement, and we
use the fractional P\'olya--Szeg\H{o} inequality
$\mathcal E_s(f^*)\leq\mathcal E_s(f)$, see \cite{AL}. Let $U^*(\cdot,z)$ be the
Schwarz rearrangement of $U(\cdot,z)$ in the horizontal variables, let
$\widehat U=\operatorname{Ext}_s(u_\Omega^*)$, and let $u_B$ and $U_B$
be the torsion function of $B$ and its extension.

For $z>0$, set $M(z)=\max_xU(x,z)$ and
\[
 \mathcal D=\{(t,z):z>0,\ 0<t<M(z)\}.
\]
For $(t,z)\in\mathcal D$ define
\[
 E_{t,z}=\{U(\cdot,z)>t\},\qquad m(t,z)=|E_{t,z}|,
 \qquad \rho(t,z)=\left(\frac{m(t,z)}{\omega_n}\right)^{1/n},
 \qquad m(t(\rho,z),z)=\omega_n\rho^n.
\]
By Lemma~\ref{lem:half-four-deficits} below, $m(\cdot,z)$ is continuous
and strictly decreasing, so that $t(\rho,z)$ is well defined for every
$\rho>0$, and $t(\rho,z)=U^*(\rho e_1,z)$. We also set
\begin{equation}\label{eq:level-notation}
 E_{\rho,z}=E_{t(\rho,z),z},\qquad
 w(\rho,z)=-\partial_\rho t(\rho,z),\qquad
 P_\rho=n\omega_n\rho^{n-1}.
\end{equation}
For the ball profile, put
\begin{equation}\label{eq:corner-notation}
 \begin{aligned}
 t_B(\rho,z)&=U_B(\rho e_1,z),&
 w_B(\rho,z)&=-\partial_\rho t_B(\rho,z),\\
 v(\rho,z)&=t(\rho,z)-t_B(\rho,z)&
 \end{aligned}
\end{equation}
Notice that 
\[
w=w_B-v_\rho.
\]
For $0<\ell<1/4$, write
\begin{equation}\label{def:R}
 \mathcal R=(1-\ell,1)\times(0,\ell).
\end{equation}

We say that $(t,z)\in\mathcal D$, or the corresponding level, is
regular if $\nabla_xU(\cdot,z)$ does not vanish on
$\{U(\cdot,z)=t\}$; a pair $(\rho,z)$ is regular if
$(t(\rho,z),z)$ is.  With $g=|\nabla_xU|$ and $h=U_z$, set
\begin{equation}\label{eq:surface-quantities}
 \begin{gathered}
 \alpha=\int_{\partial E_{t,z}}\frac1g\,d\mathcal H^{n-1},\qquad
 \beta=\int_{\partial E_{t,z}}g\,d\mathcal H^{n-1},\\
 \gamma=\int_{\partial E_{t,z}}\frac hg\,d\mathcal H^{n-1},\qquad
 \eta=\int_{\partial E_{t,z}}\frac{h^2}{g}\,d\mathcal H^{n-1}.
 \end{gathered}
\end{equation}
For the local torsion problem, Gauss--Green and $-\Delta u=1$ give
$\beta(t)=m(t)$ which is used in
\cite[Proposition~2.5]{GMR} to simplify some computations. Here we cannot rely on such a simplification.

Let $P=\operatorname{Per}(E_{t,z})$ and
$P_*=P_{\rho(t,z)}$. Define
\begin{equation}\label{eq:four-deficits}
 \begin{aligned}
 D_H&=\alpha\beta-P^2,&D_I&=P^2-P_*^2,\\
 D_V&=\alpha\eta-\gamma^2,&
 D_B&=\left|\binom{-P_*}{\gamma}
          -\alpha\nabla t_B(\rho(t,z),z)\right|^2.
 \end{aligned}
\end{equation}
The quotients involving these quantities are set equal to zero on
nonregular levels. In radius coordinates they are evaluated at
$(t(\rho,z),z)$.

Finally, for a cap $A>0$, define
\begin{equation}\label{def:Q}
 \begin{gathered}
 Q(t,z)=\inf_{a\in\R^n}|E_{t,z}\mathbin\triangle
 B_{\rho(t,z)}(a)|,\qquad
 q(\rho,z)=Q(t(\rho,z),z),\\
 Q_A=\min\{Q,A\},\qquad q_A=\min\{q,A\}.
 \end{gathered}
\end{equation}
The cap plays a role only in the proof of \eqref{eq:trace-result}, where
$A=\omega_n\mathcal A(\Omega)$: it makes the error terms there
proportional to $A$, see the discussion before \eqref{eq:unified-window}.
At regular levels, the fixed-volume direction is
\begin{equation}\label{def:Z}
 \mathcal Z=\partial_z+\frac\gamma\alpha\partial_t,
 \qquad \mathcal Zm=0.
\end{equation}

\subsection{Horizontal rearrangement of the extension}

\begin{lemma}\label{lem:rearranged-extension}
Let $\Omega\subset\R^n$ be bounded and open, with $|\Omega|=\omega_n$,
and let $U=\operatorname{Ext}_su_\Omega$. For every $z>0$, the function
$U(\cdot,z)$ is nonnegative, belongs to $H^1(\R^n)$, and vanishes at
infinity. It is positive and real analytic, its positive
superlevel sets are bounded, and its level sets and horizontal critical
set are Lebesgue null. The horizontal rearrangement $U^*$ has weak derivatives in
$\R^{n+1}_+$ and, for almost every $z>0$,
\begin{equation}\label{eq:slice-contraction}
 \|\nabla_xU^*(\cdot,z)\|_2\leq\|\nabla_xU(\cdot,z)\|_2,\qquad
 \|\partial_zU^*(\cdot,z)\|_2\leq\|U_z(\cdot,z)\|_2.
\end{equation}
In particular $U^*$ has finite weighted Dirichlet energy. Moreover, as $z\downarrow0$,
\begin{equation}\label{eq:rearranged-trace}
 \|U^*(\cdot,z)-u_\Omega^*\|_2\leq\|U(\cdot,z)-u_\Omega\|_2
 \longrightarrow0,
\end{equation}
so that the bottom trace of $U^*$ is $u_\Omega^*$.
\end{lemma}

\begin{proof}
Write $u=u_\Omega$. It is nonnegative, belongs to $L^1\cap L^2$, and has
bounded support. For $z>0$ the kernel $P_z^{(s)}$ and its gradient are
integrable, so $U(\cdot,z)=P_z^{(s)}*u$ and
$\nabla_xU(\cdot,z)=(\nabla P_z^{(s)})*u$ belong to $L^2(\R^n)$, and
$U(\cdot,z)$ tends to zero at infinity.  Real analyticity follows from the formula for $P_z^{(s)}$.
Its positive superlevel sets are bounded, uniformly when $z$ ranges in
a compact subinterval of $(0,\infty)$ and $t$ stays bounded away from zero.
Since at least one coordinate derivative is a nonzero analytic function,
\[
 |\{x:\nabla_xU(x,z)=0\}|=0.
\]
Likewise, every level set has zero $n$-dimensional measure.

Horizontal
P\'olya--Szeg\H{o} gives the first inequality in
\eqref{eq:slice-contraction}. The Hardy--Littlewood--P\'olya inequality
gives
\[
 \|U^*(\cdot,z+\varepsilon)-U^*(\cdot,z)\|_2
 \leq\|U(\cdot,z+\varepsilon)-U(\cdot,z)\|_2.
\]
On every compact height interval, $z\mapsto U^*(\cdot,z)$ is therefore
Lipschitz with values in $L^2(\R^n)$, and difference quotients give a
weak vertical derivative satisfying the second inequality in
\eqref{eq:slice-contraction} for almost every $z$. Multiplying these
slice inequalities by $z^{1-2s}$ and integrating shows that $U^*$ has
finite weighted Dirichlet energy. On compact positive height intervals
the weight is bounded above and below. By \eqref{eq:slice-contraction},
$U^*\in L^2_{\rm loc}((0,\infty);H^1(\R^n))$ and $\partial_zU^*\in L^2_{\rm loc}((0,\infty);L^2(\R^n))$,
so by Fubini the slice-wise gradient $\nabla_xU^*$ and the weak vertical derivative are the joint weak derivatives, and
$U^*\in H^1_{\rm loc}(\R^n\times(0,\infty))$. Finally, the
kernel $P_z^{(s)}$ is an approximate identity, so $U(\cdot,z)\to u$ in
$L^2(\R^n)$, and the Hardy--Littlewood--P\'olya inequality gives
\eqref{eq:rearranged-trace}.
\end{proof}

\subsection{Exact torsion and extension splittings}

\begin{lemma}\label{lem:split}
Let $\Omega\subset\R^n$ be bounded and open, with $|\Omega|=\omega_n$,
let $u=u_\Omega$ and $U=\operatorname{Ext}_su$, and let $U^*$ and $U_B$
be as above. Then
\begin{equation}\label{eq:torsion-split}
 \delta_T(\Omega)=\mathcal E_s(u)-\mathcal E_s(u^*)
                +\mathcal E_s(u^*-u_B).
\end{equation}
The difference $\mathcal E_s(u)-\mathcal E_s(u^*)$ and the last term
are nonnegative. In extension variables, with $d_s$ as in
\eqref{eq:general-extension-normalization}, one also has
\begin{equation}\label{eq:extension-split-new}
 \delta_T(\Omega)=d_s\int_{\R^{n+1}_+}z^{1-2s}
 \left(|\nabla U|^2-|\nabla U^*|^2
             +|\nabla(U^*-U_B)|^2\right)\,dx\,dz.
\end{equation}
In particular, by \eqref{eq:slice-contraction}
\begin{equation}\label{eq:split-weighted-profile-control}
 \int_{\R^{n+1}_+}z^{1-2s}|\nabla(U^*-U_B)|^2\,dx\,dz
 \leq C\delta_T(\Omega).
\end{equation}
\end{lemma}

\begin{proof}
Fractional P\'olya--Szeg\H{o} gives
$\mathcal E_s(u^*)\leq\mathcal E_s(u)$.
The rearranged function belongs to $H_0^s(B)$: it is in $H^s(\R^n)$
and is supported in $\overline B$. For $0<r<1$, the inward dilation
$u^*(x/r)$ is supported in $\overline{B_r}$ and converges to $u^*$ in
$H^s(\R^n)$ as $r\uparrow1$. Mollification at a scale smaller than
$1-r$ then gives an approximation by $C_c^\infty(B)$ functions.
Thus it is an admissible test in the ball torsion equation under the
fixed closure convention. We then have
\[
 \mathcal E_s(u_B,u^*)=\int_Bu^*=T_s(\Omega),\qquad
 \mathcal E_s(u_B)=T_s(B),\qquad \mathcal E_s(u)=T_s(\Omega).
\]
Consequently, after expanding the squares,
\[
 \mathcal E_s(u)-\mathcal E_s(u^*)+\mathcal E_s(u^*-u_B)
 =T_s(\Omega)+T_s(B)-2T_s(\Omega)=\delta_T(\Omega),
\]
which proves \eqref{eq:torsion-split} and the nonnegativity of its
 two grouped terms.

In this proof write
\[
 \langle F,G\rangle
 =d_s\int_{\R^{n+1}_+}z^{1-2s}\nabla F\cdot\nabla G\,dx\,dz,
 \qquad \|F\|^2=\langle F,F\rangle.
\]
 In particular,
$\|\operatorname{Ext}_s f\|^2=\mathcal E_s(f)$.
Recall that $\widehat U=\operatorname{Ext}_s(u^*)$.
Lemma~\ref{lem:rearranged-extension} shows that $U^*$ has finite weighted
energy and bottom trace $u^*$. Thus $U^*-\widehat U$ has zero weighted trace.
For any finite-energy weighted-harmonic extension $G$, the variation
$G+\varepsilon(U^*-\widehat U)$ has the same trace as $G$.
Using the variational characterization \eqref{varextension} of $\mathcal E_s$ and differentiation 
at $\varepsilon=0$ give
\begin{equation}\label{eq:split-zero-trace-orthogonality}
 \langle G,U^*-\widehat U\rangle=0.
\end{equation}
 Apply \eqref{eq:split-zero-trace-orthogonality} to $G=\widehat U$
and to $G=\widehat U-U_B=\operatorname{Ext}_s(u^*-u_B)$. We obtain
\begin{equation}\label{eq:split-extension-pythagoras}
 \begin{aligned}
 \|U^*\|^2&=\|\widehat U\|^2+\|U^*-\widehat U\|^2,\\
 \|U^*-U_B\|^2&=\|U^*-\widehat U\|^2+\|\widehat U-U_B\|^2.
 \end{aligned}
\end{equation}
Substitution into \eqref{eq:torsion-split} cancels the common
$\|U^*-\widehat U\|^2$ term and gives
\[
 \delta_T(\Omega)=\|U\|^2-\|U^*\|^2+\|U^*-U_B\|^2,
\]
which is \eqref{eq:extension-split-new}.
\end{proof}

\subsection{Ball speed and boundary profile}
Two further estimates are used in Section~\ref{sec:bounded-torsion}; since
their proofs are independent of the rest of the argument, they are given
in Appendix~\ref{app:ball-hardy}. Lemma~\ref{lem:half-ball-speed} shows
that there is $\ell\in(0,1/16)$, depending only on $n$
and $s$, such that the radial speed $w_B$ of the ball extension
satisfies, on the rectangle $\mathcal R$ of \eqref{def:R},
\[
 w_B(\rho,z)\asymp(1-\rho+z)^{s-1},\qquad
 |\nabla w_B(\rho,z)|\lesssim\frac{w_B(\rho,z)}{1-\rho+z}.
\]
Lemma~\ref{lem:hardy} shows that the profile error $v=t-t_B$ satisfies
\[
 \int_{\mathcal R}z^{1-2s}|\nabla v|^2\,d\rho\,dz
 +\int_{1-\ell}^1\frac{|v(\rho,0)|^2}{(1-\rho)^{2s}}\,d\rho
 \leq C(\ell)\delta_T(\Omega),
\]
where $v(\cdot,0)$ is the bottom trace of $v$.

\section{Quantitative torsion inequality for bounded sets}
\label{sec:bounded-torsion}

\subsection{Four-defect identity}

The next lemma is the counterpart of the level-set identity
\cite[Theorem~2.4]{GMR} for the local torsion function and of its
expression in radius coordinates \cite[Proposition~2.5]{GMR}. The defects
$D_H$ and $D_I$ are those of \cite{GMR}, computed on the horizontal
slices of the extension; the vertical defect $D_V$ and the profile defect
$D_B$ have no local counterpart.

Using $\alpha w=P_\rho$, proved below, the two horizontal defects
can still be combined as
\[
 \frac{w}{P_\rho^2}(D_H+D_I)=\frac{\beta}{P_\rho}-w\geq0.
\]
In the local case, $\beta=m=\omega_n\rho^n$ reduces the right-hand
side to $\rho/n-w$, the nonnegative difference from the ball speed.
For the extension, $\beta/P_\rho$ is not prescribed by the level-set
volume, as explained after \eqref{eq:surface-quantities}. This identity
therefore does not give a pointwise comparison between $w$ and $w_B$.
The profile defect $D_B$ provides a nonnegative energy control of
their difference instead.

\begin{lemma}\label{lem:half-four-deficits}
Let $\Omega\subset\R^n$ be bounded and open, with $|\Omega|=\omega_n$,
let $u=u_\Omega$ and $U=\operatorname{Ext}_su$.  The function $m(\cdot,z)$ is continuous, strictly
decreasing and locally absolutely continuous on $(0,M(z))$, with
$m(t,z)\to+\infty$ as $t\downarrow0$ and $m(t,z)\to0$ as
$t\uparrow M(z)$. The nonregular pairs form a null set, both in the
$(t,z)$ and in the $(\rho,z)$ variables, and at regular pairs
\begin{equation}\label{eq:half-level-derivatives}
 m_t=-\alpha,\qquad m_z=\gamma,\qquad
 \partial_\rho t=-\frac{P_\rho}{\alpha},\qquad
 \partial_z t=\frac{\gamma}{\alpha}.
\end{equation}
All four defects are nonnegative and
\begin{equation}\label{eq:half-four-deficit-identity}
 \delta_T(\Omega)=T_s(B)-T_s(\Omega)
 =d_s\int_0^\infty z^{1-2s}\int_0^{M(z)}
       \frac{D_H+D_I+D_B+D_V}{\alpha}\,dt\,dz.
\end{equation}
More precisely,
\begin{align}
 \int_{\R^{n+1}_+}z^{1-2s}\bigl(|\nabla_xU|^2-|\nabla_xU^*|^2\bigr)\,dx\,dz
 &=\int_0^\infty z^{1-2s}\int_0^{M(z)}\frac{D_H+D_I}{\alpha}\,dt\,dz,
 \label{eq:half-horizontal-loss}\\
 \int_{\R^{n+1}_+}z^{1-2s}\bigl(|U_z|^2-|\partial_zU^*|^2\bigr)\,dx\,dz
 &=\int_0^\infty z^{1-2s}\int_0^{M(z)}\frac{D_V}{\alpha}\,dt\,dz,
 \label{eq:half-vertical-loss}\\
 \int_{\R^{n+1}_+}z^{1-2s}|\nabla(U^*-U_B)|^2\,dx\,dz
 &=\int_0^\infty z^{1-2s}\int_0^{M(z)}\frac{D_B}{\alpha}\,dt\,dz.
 \label{eq:half-ball-loss}
\end{align}
In radius coordinates this reads
\begin{equation}\label{eq:half-four-deficit-radial}
 \delta_T(\Omega)=d_s\int_0^\infty z^{1-2s}\int_0^\infty
 \frac{w^2}{P_\rho}(D_H+D_I+D_B+D_V)\,d\rho\,dz.
\end{equation}
Consequently, on the fixed rectangle $\mathcal R$,
\begin{equation}\label{eq:half-collar-budget}
 \int_{\mathcal R}z^{1-2s}\left[
 \frac{w^2}{P_\rho}(D_H+D_I+D_V)+P_\rho|\nabla v|^2
 \right]\,d\rho\,dz\leq C\delta_T(\Omega).
\end{equation}
\end{lemma}

\begin{proof}

 Coarea gives
\begin{equation}\label{eq:half-volume-ac}
 m(a,z)-m(b,z)=\int_a^b\alpha(t,z)\,dt,
 \qquad 0<a<b<M(z).
\end{equation}
The function $m(\cdot,z)$ is locally absolutely continuous and strictly
decreasing: the preimage of every nonempty interval of values inside
$(0,M(z))$ contains a nonempty open set. Its endpoint limits are
$+\infty$ at $t=0$ and zero at $t=M(z)$. Hence $t(\rho,z)$ is uniquely
defined for every $\rho>0$.

At a regular pair $(t,z)$ the  outward normal to $\partial E_{t,z}$ is $-\nabla_xU/g$. At fixed $t$, its
outward normal velocity is $h/g$, obtained by differentiating $U(x,z)=t$.
This gives the first two identities in \eqref{eq:half-level-derivatives};
the last two follow by differentiating $m(t(\rho,z),z)=\omega_n\rho^n$.

For fixed $z$, coarea and the nullity of the horizontal critical set give
\[
 \int_{\R^n}|\nabla_xU|^2\,dx=\int_0^{M(z)}\beta\,dt,
 \qquad
 \int_{\R^n}|U_z|^2\,dx=\int_0^{M(z)}\eta\,dt.
\]
The monotone map $t\mapsto\rho(t,z)$ is locally absolutely continuous,
with $|\partial_t\rho|=\alpha/P_*$. Polar integration and
\eqref{eq:half-level-derivatives} therefore give, for a.e. $z$,
\[
 \int_{\R^n}|\nabla_xU^*|^2\,dx
      =\int_0^{M(z)}\frac{P_*^2}{\alpha}\,dt,
 \qquad
 \int_{\R^n}|\partial_zU^*|^2\,dx
      =\int_0^{M(z)}\frac{\gamma^2}{\alpha}\,dt.
\]
Multiplying by $z^{1-2s}$, integrating in height, and subtracting proves \eqref{eq:half-horizontal-loss} and
\eqref{eq:half-vertical-loss}, since
$\alpha\beta-P_*^2=D_H+D_I$ and $\alpha\eta-\gamma^2=D_V$.

For the profile error, both $U^*$ and $U_B$ are radial in $x$.
Using the same change of variables, now with
$P_\rho\,|d\rho|=\alpha\,dt$, gives
\begin{align*}
 \int_{\R^{n+1}_+}z^{1-2s}|\nabla(U^*-U_B)|^2\,dx\,dz
 &=\int_0^\infty z^{1-2s}\int_0^\infty P_\rho|\nabla v(\rho,z)|^2\,d\rho\,dz\\
 &=\int_0^\infty z^{1-2s}\int_0^{M(z)}
       \alpha|\nabla v(\rho(t,z),z)|^2\,dt\,dz\\
 &=\int_0^\infty z^{1-2s}\int_0^{M(z)}\frac{D_B}{\alpha}\,dt\,dz.
\end{align*}
The last equality follows from
$\alpha\nabla t=(-P_*,\gamma)$ at the corresponding radius. This proves
\eqref{eq:half-ball-loss}. All changes of variables can first be made on
compact level intervals and then extended by monotone convergence of
the nonnegative energy densities.

Combining \eqref{eq:half-horizontal-loss}, \eqref{eq:half-vertical-loss}
and \eqref{eq:half-ball-loss} with the splitting
\eqref{eq:extension-split-new} of Lemma~\ref{lem:split} proves
\eqref{eq:half-four-deficit-identity}. Cauchy--Schwarz gives $D_H,D_V\geq0$, classical isoperimetry gives
$D_I\geq0$, and $D_B$ is a squared norm. Finally,
$|dt|=w\,d\rho$, $\alpha w=P_\rho$, and
$D_B=\alpha^2|\nabla v|^2$ give
\eqref{eq:half-four-deficit-radial} and, by restriction to $\mathcal R$,
\eqref{eq:half-collar-budget}.
\end{proof}

\subsection{A nonsharp quantitative Saint-Venant inequality}

For the reader's convenience, we use Lemma~\ref{lem:half-four-deficits}
to give a short derivation of the nonsharp quantitative Saint-Venant
inequality of \cite{BCV,Kim}.
It uses only the isoperimetric defect and the profile defect, and
provides the initial smallness needed in
Proposition~\ref{prop:half-truncation}.

\begin{proposition}\label{prop:nonsharp-saint-venant}
For every $n\geq2$ and $s\in(0,1)$ there is $c=c(n,s)>0$ such that
every open set $\Omega\subset\R^n$ with $|\Omega|=\omega_n$ satisfies
\begin{equation}\label{eq:nonsharp-saint-venant}
 T_s(B)-T_s(\Omega)\geq c\,\mathcal A(\Omega)^{2+1/s}.
\end{equation}
\end{proposition}

\begin{proof}
By truncation we may assume that $\Omega$ is bounded (but without a universal bound on the diameter). We write
$\delta=\delta_T(\Omega)$ and $\mathcal A=\mathcal A(\Omega)$.
We may assume that $\mathcal A>0$ since otherwise there is nothing to prove. Fix
$\ell$ of Lemma~\ref{lem:half-ball-speed}, and put
\[
 \tau=\eps\mathcal A,\qquad
 I_\tau=(1-2\tau,1-\tau),\qquad
 t_\tau=b\tau^s,\qquad z_\tau=\tau^{1+1/(2s)}.
\]
Here $b=b(n,s)>0$ is chosen small enough below, and then
$\eps=\eps(n,s)>0$ is chosen so that $2\tau<\ell$ and the lower bound
\eqref{eq:nonsharp-level-asymmetry} below holds. Since $\mathcal A\leq2$, we may also assume
$\tau<1$; hence $z_\tau\leq\tau$.

Let $G=\{u>2t_\tau\}$. If $|G|<\omega_n(1-\tau)^n$, then
$u^*(\rho e_1)\leq2t_\tau$ for almost every $\rho\geq1-\tau$.
The ball formula \eqref{eq:general-ball-torsion} gives
$u_B(\rho e_1)\geq c\tau^s$ on
$(1-\tau,1-\tau/2)$. Choose $b$ so that $2t_\tau$ is at most half
this lower bound. Lemma~\ref{lem:hardy} then gives
\[
 \delta\geq c\int_{1-\tau}^{1-\tau/2}
 \frac{|u^*(\rho e_1)-u_B(\rho e_1)|^2}{(1-\rho)^{2s}}\,d\rho
 \geq c\tau\geq c\mathcal A^{2+1/s}.
\]
We may therefore assume that
\begin{equation}\label{eq:nonsharp-low-level-mass}
 |\Omega\setminus G|\leq\omega_n\bigl(1-(1-\tau)^n\bigr)\les \tau.
\end{equation}

The weighted Cauchy--Schwarz inequality along vertical lines and the
strong $L^2$ trace give for fixed $z>0$,
\begin{equation}\label{eq:nonsharp-extension-trace}
 \begin{aligned}
 \|U(\cdot,z)-u\|_2^2
 &\leq C z^{2s}\int_0^z\zeta^{1-2s}
                         \|U_z(\cdot,\zeta)\|_2^2\,d\zeta\\
 &\les  z^{2s}\mathcal E_s(u)\le  z^{2s}T_s(B).
 \end{aligned}
\end{equation}
For $\rho\in I_\tau$ and $0<z<z_\tau$, let $E=E_{\rho,z}$ and
$t=t(\rho,z)$. If $t\leq t_\tau$, then
$|G\setminus E|\leq t_\tau^{-2}\|U(\cdot,z)-u\|_2^2$.
Since $|E|=\omega_n\rho^n\leq|\Omega|$, we have
$|E\setminus\Omega|\leq|\Omega\setminus E|$; using
\eqref{eq:nonsharp-low-level-mass}, we obtain
\[
 |\Omega\triangle E|\les \tau+
                  t_\tau^{-2}\|U(\cdot,z)-u\|_2^2.
\]
If $t>t_\tau$, since $u=0$ and  $U(\cdot,z)>t>t_\tau$ on $E\setminus\Omega$, Chebyshev yields
$|E\setminus\Omega|\leq t_\tau^{-2}\|U(\cdot,z)-u\|_2^2$.
Now
$|\Omega\triangle E|=|\Omega|-|E|+2|E\setminus\Omega|$ and
$|\Omega|-|E|=\omega_n(1-\rho^n)\les \tau$ on $I_\tau$ give the
same bound. Thus, in both cases,
\begin{equation}\label{eq:nonsharp-level-closeness}
 |\Omega\triangle E_{\rho,z}|
 \les\tau+\frac{z^{2s}}{\tau^{2s}}\les \tau.
\end{equation}
For every center $a\in\R^n$, the triangle inequality yields
\[
 |E_{\rho,z}\triangle B_\rho(a)|
 \geq|\Omega\triangle B_1(a)|-|\Omega\triangle E_{\rho,z}|
                -\omega_n(1-\rho^n).
\]
Taking the infimum over $a$ and then choosing $\eps$ small enough gives
\begin{equation}\label{eq:nonsharp-level-asymmetry}
 q(\rho,z)\geq\omega_n\mathcal A-C\tau
             \geq\tfrac12\omega_n\mathcal A
 \qquad\text{on }I_\tau\times(0,z_\tau).
\end{equation}

The quantitative isoperimetric inequality  gives
$D_I\ges \rho^{2n-2}\mathcal A(E_{\rho,z})^2$.
Since $q=\omega_n\rho^n\mathcal A(E_{\rho,z})$, it follows that
$q^2\les \rho^2D_I$. Hence
\eqref{eq:nonsharp-level-asymmetry} gives $D_I\ges \mathcal A^2$
on the rectangle above.
Retaining $D_I$ and the radial part of $D_B$ in
\eqref{eq:half-collar-budget}, and using $P_\rho\asymp1$ there, we find
\[
 \delta\ges \int_0^{z_\tau}z^{1-2s}\int_{I_\tau}
          \bigl(\mathcal A^2w^2+(w_B-w)^2\bigr)\,d\rho\,dz.
\]
For $0<\mathcal A\leq2$,
$\mathcal A^2w^2+(w_B-w)^2\ges \mathcal A^2w_B^2$.
Lemma~\ref{lem:half-ball-speed} gives $w_B\ges \tau^{s-1}$ on
this rectangle. Consequently,
\[
 \delta\ges \mathcal A^2\tau^{2s-1}z_\tau^{2-2s}
       =\mathcal A^2\tau^{1/s}
       \ges \mathcal A^{2+1/s}.
\]

\end{proof}

\subsection{Horizontal estimates}

Throughout this subsection $\Omega\subset\R^n$ is bounded and open,
with $|\Omega|=\omega_n$, $u=u_\Omega$, $U=\operatorname{Ext}_su$, and
$E_{\rho,z}$, $q(\rho,z)$ are as in Section~\ref{sec:preliminaries}.

In \cite{GMR} the level sets of the torsion function are contained in
$\Omega$, hence in $B_R$ when $\Omega\subset B_R$; this confinement enters
the strong isoperimetric estimate \cite[Proposition~3.5]{GMR}. At positive
height $U(\cdot,z)$ is positive on all of $\R^n$, and the confinement of
the level sets of volume at most $\omega_n$ has to be proved.

\begin{lemma}\label{lem:confinement}
Let $\Omega\subset B_R$ be open with $|\Omega|=\omega_n$, and let $z_0>0$. Set $L=R+2\sqrt{4R^2+z_0^2}$. Then 
\begin{equation}\label{eq:general-level-confinement}
 E_{\rho,z}\subset B_L,
 \qquad 0<\rho\leq1,\quad 0<z\leq z_0.
\end{equation}
\end{lemma}

\begin{proof}
Since $|\Omega|=\omega_n$ and $\Omega\subset B_R$, necessarily $R\geq1$.
The torsion function is nonnegative, supported in $\overline\Omega$,
and $T_s(\Omega)=\int u>0$. By \eqref{eq:general-poisson}, for
$x,y\in B_R$ and $0<z\leq z_0$,
\[
 P_z^{(s)}(x-y)
 \geq\frac{\Gamma(n/2+s)}{\pi^{n/2}\Gamma(s)}
       \frac{z^{2s}}{(4R^2+z_0^2)^{(n+2s)/2}}.
\]
Integrating against $u(y)$ therefore gives, throughout $B_R$,
\begin{equation}\label{eq:confinement-interior-lower}
 U(x,z)\geq
 \frac{\Gamma(n/2+s)}{\pi^{n/2}\Gamma(s)}
 \frac{z^{2s}T_s(\Omega)}{(4R^2+z_0^2)^{(n+2s)/2}}=t_z.
\end{equation}
Since $B_R\subset E_{t_z,z}$, we have 
\[
|E_{t_z,z}|\ge |B_R|\ge \omega_n=|E_{t(1,z),z}|
\] 
and thus 
\begin{equation}\label{eq:confinement-threshold-lower}
 t(1,z)\geq t_z.
\end{equation}

For $|x|\geq L>R$ and $y\in\Omega$, we have $|x-y|\geq L-R$.
Dropping the nonnegative $z^2$ from the kernel denominator yields
\begin{equation}\label{eq:confinement-exterior-upper}
 U(x,z)\leq
 \frac{\Gamma(n/2+s)}{\pi^{n/2}\Gamma(s)}
 \frac{z^{2s}T_s(\Omega)}{(L-R)^{n+2s}}.
\end{equation}
For  $L=R+2\sqrt{4R^2+z_0^2}$, the right-hand side of
\eqref{eq:confinement-exterior-upper} is then $2^{-(n+2s)}$ times
the right-hand side of \eqref{eq:confinement-threshold-lower}.
Since $2^{-(n+2s)}<1$, it follows that $U(x,z)<t(1,z)$ whenever
$|x|\geq L$. This confines $E_{1,z}$ to $B_L$ for every
$0<z\leq z_0$. Finally, $t(\rho,z)\geq t(1,z)$ for $\rho\leq1$,
so $E_{\rho,z}\subset E_{1,z}$ and
\eqref{eq:general-level-confinement} follows.
\end{proof}
We now record the horizontal estimates. At each fixed height, they follow from items (1)--(3) of
\cite[Proposition~3.6]{GMR} together with the estimate on good radii in the
proof of \cite[Theorem~3.7]{GMR}. It builds on \cite[Lemmas~3.1 and~3.2]{GMR}.
The only additional input  in our setting is the
confinement from Lemma~\ref{lem:confinement}. Notice that smoothness of $\partial\Omega$ is assumed in
\cite[Lemma~3.2]{GMR} only to obtain a uniform first-order expansion of
the torsion function near the level sets. Here $U(\cdot,z)$ is smooth on
$\R^n$, so the argument carries over. We refer the reader to \cite{GMR} for the proofs.

\begin{proposition}
\label{prop:half-horizontal-moving-level}
Let $\Omega\subset B_R$  with $R\geq1$. All defects below are evaluated
at $(t(\rho,z),z)$. For every $z>0$ the infimum defining
$q(\rho,z)$ is attained for every
$\rho>0$, and for
$0<\rho,\rho'\leq1$,
\begin{equation}\label{eq:qlip}
 |q(\rho,z)-q(\rho',z)|
 \leq2\omega_n|\rho^n-(\rho')^n|
 \leq2n\omega_n|\rho-\rho'|,
 \qquad 0\leq q(\rho,z)\leq2\omega_n.
\end{equation}
In particular, $q(\cdot,z)$ is locally Lipschitz on $(0,\infty)$.

Fix  $Z>0$, and $\rho_*\in(0,1)$.  There is $C=C(R,Z,\rho_*)>0$ such that
\begin{equation}
\label{eq:half-horizontal-pointwise-action}
 q(\rho,z)^2+\bigl(\partial_\rho q(\rho,z)\bigr)_+^2\leq C(D_H+D_I)
\end{equation}
for almost every $(\rho,z)\in[\rho_*,1]\times(0,Z)$. Consequently,
\begin{equation}
\label{eq:half-horizontal-action-budget-unweighted}
 \int_0^Z z^{1-2s}\int_{\rho_*}^1
 w^2\left(q^2+(\partial_\rho q)_+^2\right)\,d\rho\,dz
 \leq C\delta_T(\Omega).
\end{equation}
\end{proposition}

\subsection{Scalar differentiation in level coordinates}

We now turn to the treatment of the vertical derivatives, which is new compared with \cite{GMR}. For this, we first revert to the original $(t,z)$ variables.

\begin{lemma}\label{lem:half-level-ACL}
Let $\Omega\subset\R^n$ be bounded and open, with $|\Omega|=\omega_n$,
let $u=u_\Omega$ and $U=\operatorname{Ext}_su$, and let $A>0$. The
functions $m$, $Q$, $Q_A$ of \eqref{def:Q} are continuous on $\mathcal D$, and they belong to
$W^{1,1}_{\rm loc}(\mathcal D)$.
Writing
\[
 \vartheta(t,z)=\int_{\partial E_{t,z}}\frac{|h|}{g}
 \,d\mathcal H^{n-1}
\]
we have 
\begin{equation}\label{eq:half-level-weak-bounds}
 |(Q_A)_t|\leq2\alpha,\quad |(Q_A)_z|\leq2\vartheta.
\end{equation}
Furthermore, recalling the definition \eqref{def:Z} of $\mathcal Z$,
\begin{equation}\label{eq:half-fixed-volume-direction}
 |\mathcal ZQ_A|^2\leq D_V\qquad\text{a.e. in }\mathcal D.
\end{equation}
\end{lemma}

\begin{proof}

\emph{Step 1: the asymmetry is $2$-Lipschitz for the symmetric
difference.} If $E,F$ have finite measure and $B_E,B_F$ denote balls
with $|B_E|=|E|$, $|B_F|=|F|$, comparison with a common center gives
\begin{equation}\label{eq:Q-lipschitz}
 \left|\inf_a|E\triangle B_E(a)|-
             \inf_a|F\triangle B_F(a)|\right|
 \leq |E\triangle F|+\bigl||E|-|F|\bigr|
 \leq2|E\triangle F|.
\end{equation}
Taking a cap preserves this Lipschitz bound, since $r\mapsto\min\{r,A\}$
is $1$-Lipschitz and the cap is independent of the pair. Everything in
the first two assertions follows by combining \eqref{eq:Q-lipschitz}
with estimates for $|E_{t_2,z_2}\triangle E_{t_1,z_1}|$.

\emph{Step 2: continuity.} On the rectangles under
consideration all the sets $E_{t,z}$ lie in one fixed bounded spatial
set, $U$ is continuous, and its positive level sets have zero  Lebesgue measure.
Dominated convergence therefore makes
$(t,z)\mapsto\mathbf1_{E_{t,z}}$ continuous from $\mathcal D$ to
$L^1(\R^n)$. Hence $m$ is continuous, and by \eqref{eq:Q-lipschitz} so
are $Q$ and $Q_A$. 

\emph{Step 3: horizontal absolute continuity.} At fixed $z$ the sets
$E_{t,z}$ are nested, so for $t_1<t_2$
\[
 |E_{t_1,z}\triangle E_{t_2,z}|
 =m(t_1,z)-m(t_2,z)
 =\int_{t_1}^{t_2}\alpha(t,z)\,dt
\]
by \eqref{eq:half-volume-ac}. With \eqref{eq:Q-lipschitz} this gives
absolute continuity in $t$ and $|(Q_A)_t|\leq2\alpha$.

\emph{Step 4: vertical slicing bound.} We claim that for almost every
$t$ and every $z_1<z_2$ in the rectangle,
\begin{equation}\label{eq:vertical-slicing}
 |E_{t,z_2}\triangle E_{t,z_1}|
 \leq\int_{z_1}^{z_2}\vartheta(t,z)\,dz.
\end{equation}
 For almost every $t$, let $\mu_t$ be the positive measure
\[
 \mu_t(I)=\int_{\R^n}
 |D_z\mathbf1_{\{U(x,\cdot)>t\}}|(I)\,dx.
\]
It is locally finite on the parameter rectangles under consideration.
For a nonnegative continuous compactly supported function $f(t,z)$,
one-dimensional coarea followed by horizontal coarea gives
\[
 \begin{aligned}
 \int dt\int f(t,z)\,d\mu_t(z)
 &=\int\!\!\int_{\R^n}|U_z(x,z)|f(U(x,z),z)\,dx\,dz\\
 &=\int\!\!\int f(t,z)\vartheta(t,z)\,dt\,dz.
 \end{aligned}
\]
 Disintegration of these positive
measures gives $d\mu_t(z)=\vartheta(t,z)\,dz$ for almost every $t$,
which gives \eqref{eq:vertical-slicing}.
Initially this holds at almost every pair of endpoints. It extends to
all endpoints on the interval by the $L^1_x$ continuity of Step 2.
The Lipschitz comparison \eqref{eq:Q-lipschitz} now gives
absolute continuity of $Q_A$ in $z$ and its second bound in
\eqref{eq:half-level-weak-bounds}. 

\emph{Step 5: the fixed-volume direction.} We finally prove
\eqref{eq:half-fixed-volume-direction}. Since $g$ is bounded below near
every regular pair, the bounds \eqref{eq:half-level-weak-bounds} show
that $Q_A$ is Lipschitz in a neighborhood of every regular pair; in
particular, its weak and pointwise derivatives agree almost everywhere.
Fix a regular pair $(t,z)$ where $Q_A$ is differentiable, set
$c_*=\gamma(t,z)/\alpha(t,z)$ and
$(t_\varepsilon,z_\varepsilon)=(t+\varepsilon c_*,z+\varepsilon)$, so that
\[
 \mathcal ZQ_A(t,z)=\lim_{\varepsilon\to0}
 \frac{Q_A(t_\varepsilon,z_\varepsilon)-Q_A(t,z)}{\varepsilon}.
\]
By \eqref{eq:Q-lipschitz},
\[
 |Q_A(t_\varepsilon,z_\varepsilon)-Q_A(t,z)|
 \leq|E_{t_\varepsilon,z_\varepsilon}\triangle E_{t,z}|
 +|m(t_\varepsilon,z_\varepsilon)-m(t,z)|.
\]
Since $\mathcal Zm=0$ by \eqref{def:Z},
$m(t_\varepsilon,z_\varepsilon)=m(t,z)+o(|\varepsilon|)$, so the second
term does not contribute in the limit.

We compute the first term by a normal graph. The surface
$\Sigma=\partial E_{t,z}$ is compact and of class $C^1$, with $g>0$ on
it, and its outer normal is $\nu=-\nabla_xU/g$. Parametrize a tubular
neighborhood by $x=y+\sigma\nu(y)$, $y\in\Sigma$. Since
\[
 \partial_\sigma U\bigl(y+\sigma\nu(y),z\bigr)\Big|_{\sigma=0}
 =\nabla_xU\cdot\nu=-g<0,
\]
the implicit function theorem provides, for each $y\in\Sigma$ and all
small $\varepsilon$, a unique small $r_\varepsilon(y)$ with
$U\bigl(y+r_\varepsilon(y)\nu(y),z_\varepsilon\bigr)=t_\varepsilon$.
Differentiating this identity at $\varepsilon=0$ gives
$-g\,\partial_\varepsilon r_\varepsilon|_{\varepsilon=0}+h=c_*$, that is,
\[
 r_\varepsilon(y)=\varepsilon\,\frac{h(y)-c_*}{g(y)}+o(|\varepsilon|)
 \qquad\text{uniformly on }\Sigma .
\]
Thus $E_{t_\varepsilon,z_\varepsilon}$ is, near $\Sigma$, the normal
graph of $r_\varepsilon$ over $E_{t,z}$, and the area formula in normal
coordinates gives
\[
 \lim_{\varepsilon\to0}
 \frac{|E_{t_\varepsilon,z_\varepsilon}\triangle E_{t,z}|}
 {|\varepsilon|}
 =\int_{\partial E_{t,z}}\frac{|h-c_*|}{g}\,d\mathcal H^{n-1}.
\]
Recalling that $c_*=\gamma/\alpha$, we obtain
\[
|\mathcal Z Q_A (t,z)|\le \int_{\partial E_{t,z}}\frac{|h-\gamma/\alpha|}{g}\,d\mathcal H^{n-1}.
\]
Cauchy--Schwarz concludes the proof of
\eqref{eq:half-fixed-volume-direction}.
\end{proof}

\subsection{Weighted vertical flux and weak change of variables}

\begin{proposition}\label{prop:half-vertical-flux}
Let $\ell$ be as in Lemma~\ref{lem:half-ball-speed},
$\Omega\subset\R^n$ be bounded and open with $|\Omega|=\omega_n$, and let
$\mathcal R$ be the rectangle \eqref{def:R}. For each $A\in(0,2\omega_n]$, define almost everywhere on
$\mathcal R$
\[
 J_A(\rho,z)=w(\rho,z)
    (\mathcal ZQ_A)(t(\rho,z),z).
\]
Then
\begin{equation}\label{eq:flux-budget}
 |J_A|^2\leq w^2D_V,\qquad
 \int_{\mathcal R}z^{1-2s}\frac{|J_A|^2}{P_\rho}
 \,d\rho\,dz\leq C\delta_T(\Omega).
\end{equation}
For every $\varphi\in C_c^\infty(\mathcal R)$ one has
\begin{equation}\label{eq:half-flux-pullback}
 \int_{\mathcal R}J_A\varphi
 =-\int_{\mathcal R}wq_A\varphi_z
  -\int_{\mathcal R}t_z\partial_\rho(q_A\varphi).
\end{equation}
In particular, the following identity holds in distributions:
\begin{equation}\label{eq:weak-flux}
 \begin{aligned}
 \partial_z\left(\frac{wq_A}{w_B}\right)
 ={}&\frac{J_A}{w_B}
 -\partial_\rho\left(\frac{v_zq_A}{w_B}\right)
 +v_z\partial_\rho\left(\frac{q_A}{w_B}\right)
 +\frac{v_\rho q_A\partial_z w_B}{w_B^2}.
 \end{aligned}
\end{equation}
\end{proposition}

\begin{proof}
The pointwise estimate in \eqref{eq:flux-budget} follows from
\eqref{eq:half-fixed-volume-direction}. Using the definition of  $D_V$ in
\eqref{eq:half-four-deficit-radial} and integrating gives the second bound.

A null set in $(t,z)$ pulls back to a null set in $(\rho,z)$: for almost
every fixed height, $t\mapsto\rho(t,z)$ is locally absolutely continuous,
so it maps the corresponding null section to a null set, and Fubini
applies. Thus the definition of $J_A$ is independent of the
almost-everywhere representatives of the weak derivatives of $Q_A$.
The weighted energy estimate also shows that $J_A\in L^2_{\rm loc}
(\mathcal R)$, because the height weight is bounded above and below on
compact subsets of $\mathcal R$.

We now prove \eqref{eq:half-flux-pullback}. The chain rule for Sobolev
functions and the fact that $m$ stays positive on compact subsets of $\mathcal D$ give
\[
 \rho\in W^{1,1}_{\rm loc}(\mathcal D),\qquad
 \rho_t=-\frac\alpha{P_\rho},\qquad
 \rho_z=\frac\gamma{P_\rho}.
\]
For a fixed test function $\varphi$, put
$\psi(t,z)=\varphi(\rho(t,z),z)$. It is bounded and belongs to
$W^{1,1}_{\rm loc}(\mathcal D)$.  Its weak derivatives are
\[
 \psi_t=-\frac\alpha{P_\rho}\varphi_\rho,\qquad
 \psi_z=\varphi_z+\frac\gamma{P_\rho}\varphi_\rho,
\]
where the derivatives of $\varphi$ are evaluated at $(\rho(t,z),z)$.
Since $Q_A$ and $\psi$ are bounded Sobolev functions on a common compact
subset of $\mathcal D$, their product belongs to $W^{1,1}$ there, and
the weak Leibniz rule in the $z$ variable gives
\[
 \int_{\mathcal D}(\mathcal ZQ_A)\psi
 =\int_{\mathcal D}\frac\gamma\alpha (Q_A)_t\psi
 -\int_{\mathcal D}Q_A
 \left(\varphi_z+\frac\gamma{P_\rho}\varphi_\rho\right),
\]
where the functions involving $\varphi$ are evaluated at
$(\rho(t,z),z)$.
Every term is integrable on this compact support. In particular,
$|\gamma(Q_A)_t/\alpha|\leq2|\gamma|$ by
\eqref{eq:half-level-weak-bounds}.

For each fixed $z$, make the change of variables
$\rho=\rho(t,z)$. By \eqref{eq:corner-notation}, $|dt|=w\,d\rho$, and
$q_A(\rho,z)=Q_A(t(\rho,z),z)$. Hence
\[
 (q_A)_\rho=-w\,(Q_A)_t(t(\rho,z),z)\qquad\text{a.e.}
\]
 The preceding identity consequently
becomes
\[
 \int_{\mathcal R}J_A\varphi
 =-\int_{\mathcal R}\frac\gamma\alpha(q_A)_\rho\varphi
  -\int_{\mathcal R}wq_A
       \left(\varphi_z+\frac\gamma{P_\rho}\varphi_\rho\right).
\]
Using
$\gamma w/P_\rho=\gamma/\alpha=t_z$ gives precisely
\eqref{eq:half-flux-pullback}. 

In distributions, \eqref{eq:half-flux-pullback} reads
\[
 \partial_z(wq_A)=J_A-\partial_\rho(t_zq_A)+t_z(q_A)_\rho.
\]
Insert $t_z=(t_B)_z+v_z$ and
$\partial_\rho(t_B)_z=-\partial_z w_B$, and multiply by the smooth
positive function $1/w_B$ on compact subsets of $\mathcal R$. The
product rule and $w_B-w=v_\rho$ give \eqref{eq:weak-flux}.
\end{proof}

\begin{remark}\label{rem:vertic-deriv}
Although no weak vertical derivative of $q_A$ is available,
\eqref{eq:weak-flux} gives a distributional equation for $wq_A/w_B$. Notice that \eqref{eq:weak-flux} can be rewritten as

\[\partial_z(wq_A/w_B)=-\partial_\rho B_A+C_A,\]
 where
\begin{equation}\label{def:BA-CA}
 B_A:=\frac{v_z q_A}{w_B},\qquad
 C_A:=\frac{J_A}{w_B}
+v_z\,\partial_\rho\!\left(\frac{q_A}{w_B}\right)
 +\frac{v_\rho q_A(w_B)_z}{w_B^2}.
\end{equation}
These quantities are locally integrable on $\mathcal R$.
\end{remark}

\begin{corollary}[Horizontal and vertical flux estimate]
\label{cor:half-full-flux-action}
Let $R\geq1$, let $\Omega\subset B_R$ be open with $|\Omega|=\omega_n$,
and let $\ell$ be as in Lemma~\ref{lem:half-ball-speed}. Then, for every $A\in(0,2\omega_n]$, almost everywhere
on $\mathcal R$,
\begin{equation}\label{eq:half-full-pointwise-action}
 q^2+(q_\rho)_+^2+
 |(\mathcal ZQ_A)(t(\rho,z),z)|^2
 \leq C(R,\ell)(D_H+D_I)+D_V.
\end{equation}
Consequently,
\begin{equation}\label{eq:required-action}
 \begin{split}
 &\int_{\mathcal R}z^{1-2s}
 \left[w^2\bigl(q^2+(q_\rho)_+^2\bigr)+|J_A|^2+|\nabla v|^2\right]
 \,d\rho\,dz+\int_{1-\ell}^1
       \frac{|v(\rho,0)|^2}{(1-\rho)^{2s}}\,d\rho
 \leq C(R)\delta_T(\Omega).
 \end{split}
\end{equation}
\end{corollary}

\begin{proof}
Estimate \eqref{eq:half-full-pointwise-action} follows from
\eqref{eq:half-horizontal-pointwise-action} in
Proposition~\ref{prop:half-horizontal-moving-level}, with
$\rho_*=1-\ell$ and $Z=\ell$, and from
\eqref{eq:half-fixed-volume-direction}.

Multiply \eqref{eq:half-full-pointwise-action} by
$z^{1-2s}w^2/P_\rho$ and integrate over $\mathcal R$. Since
$J_A=w(\mathcal ZQ_A)(t(\rho,z),z)$, the radial four-defect identity
\eqref{eq:half-four-deficit-radial} gives
\[
 \int_{\mathcal R}z^{1-2s}\frac1{P_\rho}
 \left[w^2\bigl(q^2+(q_\rho)_+^2\bigr)+|J_A|^2\right]
 \,d\rho\,dz
 \leq C(R,\ell)\delta_T(\Omega).
\]
On $1-\ell<\rho<1$, the radial factor
$P_\rho=n\omega_n\rho^{n-1}$ is bounded above and below by positive
constants depending only on $n$ and $\ell$. Thus the preceding estimate
controls the first two bulk terms of \eqref{eq:required-action} without
the factor $P_\rho^{-1}$. Finally, Lemma~\ref{lem:hardy} gives
\[
 \int_{\mathcal R}z^{1-2s}|\nabla v|^2\,d\rho\,dz
 +\int_{1-\ell}^1
       \frac{|v(\rho,0)|^2}{(1-\rho)^{2s}}\,d\rho
 \leq C(\ell)\delta_T(\Omega).
\]
This proves \eqref{eq:required-action}. Since the choice of $\ell$
in Lemma~\ref{lem:half-ball-speed} depends only on $n$ and $s$, all
$\ell$-dependence is absorbed into $C(R)$.
\end{proof}

\subsection{The corner value of the moving-ball distance}

The asymmetry of $\Omega$ is recovered from $q$ only in the limit
$\rho\to1$, $z\to0$; this replaces the identity
$q(1)=\omega_n\mathcal A(\Omega)$ in the proof of \cite[Theorem~3.7]{GMR}.

\begin{lemma}\label{lem:corner}
Let $\Omega\subset\R^n$ be bounded and open with $|\Omega|=\omega_n$, and let
$E_z=E_{1,z}$. Then $|E_z\triangle\Omega|\to0$ as $z\downarrow0$, and
\begin{equation}\label{eq:endpoint-joint-control}
 |q(\rho,z)-\omega_n\mathcal A(\Omega)|
 \leq2n\omega_n(1-\rho)+|E_z\triangle\Omega|
 \qquad(0<\rho\leq1,\ z>0).
\end{equation}
In particular,
\begin{equation}\label{eq:qcorner}
 \lim_{\substack{\rho\to1^-\\z\to0^+}}q(\rho,z)
 =\omega_n\mathcal A(\Omega).
\end{equation}
\end{lemma}

\begin{proof}
By Lemma~\ref{lem:half-four-deficits}, $|E_z|=\omega_n=|\Omega|$ for every
$z>0$, so $|E_z\triangle\Omega|=2|\Omega\setminus E_z|$. Let $t_z=t(1,z)$.
Since $U(\cdot,z)\geq t_z$ on $E_z$ and $U(\cdot,z)\leq t_z$ outside it, and
$|\Omega|=|E_z|$, we have $\int_{E_z}U(\cdot,z)\geq\int_\Omega U(\cdot,z)$,
that is,
$\int_{\Omega\setminus E_z}U(\cdot,z)\leq\int_{E_z\setminus\Omega}U(\cdot,z)$.
Using $u=0$ outside $\Omega$, we obtain
\[
 \int_{\Omega\setminus E_z}u
 \leq\int_{\Omega\setminus E_z}|U(\cdot,z)-u|
       +\int_{E_z\setminus\Omega}|U(\cdot,z)-u|
 \leq(2\omega_n)^{1/2}\|U(\cdot,z)-u\|_2=:\eps_z,
\]
and $\eps_z\to0$ by Lemma~\ref{lem:rearranged-extension}.

We also need that $u>0$ almost everywhere in $\Omega$. Let $B'\Subset\Omega$ be a
ball. The function $f_+=(u_{B'}-u)_+$ lies in $H^s(\R^n)$ and vanishes outside $B'$, so
$f_+\in H_0^s(B')\subset H_0^s(\Omega)$. Testing the two torsion equations with $f_+$ gives
$\mathcal E_s(f,f_+)=\int_{B'}f_+-\int_\Omega f_+\leq0$, while
$\mathcal E_s(f,f_+)\geq\mathcal E_s(f_+)$. Hence $f_+=0$, that is,
$u\geq u_{B'}>0$ in $B'$. Covering $\Omega$ by countably many such balls
proves the claim.

For $\sigma>0$ we have
$\Omega\setminus E_z\subset(\Omega\cap\{u\leq\sigma\})\cup(\{u>\sigma\}\setminus E_z)$, so that
\[
 |\Omega\setminus E_z|\leq|\Omega\cap\{0<u\leq\sigma\}|+\sigma^{-1}\eps_z.
\]
Letting first $z\downarrow0$ and then $\sigma\downarrow0$ (using $|\Omega|<\infty$) gives
$|\Omega\setminus E_z|\to0$.

For every center $a$, $\bigl||E_z\triangle B_1(a)|-|\Omega\triangle B_1(a)|\bigr|\leq|E_z\triangle\Omega|$.
Taking the infimum over $a$ and recalling $|\Omega|=\omega_n$ gives
$|q(1,z)-\omega_n\mathcal A(\Omega)|\leq|E_z\triangle\Omega|$.
Together with \eqref{eq:qlip}, this proves \eqref{eq:endpoint-joint-control} and
hence \eqref{eq:qcorner}.
\end{proof}

\subsection{Quantitative torsion inequality for confined sets}

We now combine the weighted action estimate \eqref{eq:required-action} with the weak flux identity \eqref{eq:weak-flux}
 to prove the quantitative estimate for confined sets. The
following theorem is the counterpart of \cite[Theorem~3.7]{GMR}. The proof is substantially more technical than its local counterpart.

\begin{theorem}\label{thm:half-confined-SV}
For every $n\geq2$, $s\in(0,1)$, and $R\geq1$ there is
$C(R)>0$ such that every open set $\Omega\subset B_R$ with
$|\Omega|=\omega_n$ satisfies
\begin{equation}\label{eq:half-confined-SV}
 \mathcal A(\Omega)^2\leq C(R)\delta_T(\Omega).
\end{equation}
\end{theorem}

\begin{proof}
Let $\ell$, $c_0$ and $C_0$ be as in Lemma~\ref{lem:half-ball-speed}; they
depend only on $n$ and $s$.
Proposition~\ref{prop:half-horizontal-moving-level} yields that $q$ has horizontal
Lipschitz constant $\Lambda=2n\omega_n$ and $0\leq q\leq \Lambda$.
Set $A=\omega_n\mathcal A(\Omega)$; by Lemma~\ref{lem:corner},
$q(1,z)\to A$ as $z\downarrow0$. If $A=0$, the conclusion follows from $\delta_T(\Omega)\geq0$;
assume henceforth that $A>0$. Define
\begin{equation}\label{eq:J}
 \begin{split}
 \mathscr J={}&\int_{\mathcal R}z^{1-2s}
 \left[w^2\bigl(q^2+(q_\rho)_+^2\bigr)+|J_A|^2+|\nabla v|^2\right]
 \,d\rho\,dz
 +\int_{1-\ell}^1\frac{|v(\rho,0)|^2}{(1-\rho)^{2s}}\,d\rho.
 \end{split}
\end{equation}
We claim that
\begin{equation}\label{eq:trace-result}
 A^2\leq C\mathscr J,
\end{equation}
where $C$ depends only on $s,\ell,c_0,C_0$. By
Corollary~\ref{cor:half-full-flux-action},
this would conclude the proof of \eqref{eq:half-confined-SV}.

We thus prove \eqref{eq:trace-result}. Choose a nonnegative
$\chi\in C_c^\infty((1,2))$ with $\int\chi=1$, put $r_0=\ell/4$, and
fix $K\geq\max\{1,2\Lambda/r_0\}$. Write $\tau=A/K$, so
$0<\tau\leq r_0/2$. 

Let us first give an idea of the proof. Formally, one would like to compare $A$ with an interior value by
integrating a derivative in $z$. However, $q_A$ has no useful Sobolev
regularity in $z$, and bulk $L^2$ estimates do not control a prescribed
horizontal slice. The weak flux identity \eqref{eq:weak-flux} allows us
instead to differentiate horizontal averages of $wq_A/w_B$ at positive
heights. We average on windows of width $\lambda(z)=z+\tau$ near $\rho=1$. On such a
window $1-\rho+z\asymp\lambda(z)$, which matches the ball-speed estimate;
the scale regularizes the corner at size $\tau$ and becomes comparable to $z$ for
$z\gg\tau$. The corner limit and the Hardy boundary term bound the
liminf of the average at the bottom from below, a one-sided derivative estimate
propagates this bound to a fixed height interval, and there the term $w^2q^2$ in
$\mathscr J$ controls the average. The cap in $q_A$ is what makes this
work at the sharp scale. At positive heights $q$ need not be of order
$A$, while in the error terms it appears undifferentiated against
derivatives of size $\lambda^{-1}$. The bound $q_A\leq A$ turns these
factors into $A/\lambda\leq A/\tau=K$. With only $q\leq2\omega_n$ they
would be of size $\tau^{-1}$, and the argument would give
$A^3\lesssim\mathscr J$ instead of $A^2\lesssim\mathscr J$. 

We now start the proof of \eqref{eq:trace-result}. Set
\begin{equation}\label{eq:unified-window}
 \begin{gathered}
 \lambda(z)=z+\tau,\qquad
 I_{\lambda(z)}=(1-2\lambda(z),1-\lambda(z)),\\
 \varphi(\rho,z)=\frac1{\lambda(z)}
       \chi\!\left(\frac{1-\rho}{\lambda(z)}\right),\qquad
 H(z)=\int_\rho\frac{wq_A}{w_B}\varphi\,d\rho.
 \end{gathered}
\end{equation}
For $0\leq z\leq r_0$, we have $\tau\leq\lambda(z)\leq3r_0/2$,
and hence $2\lambda(z)\leq3\ell/4<\ell$. Thus the horizontal
support remains in a compact subset of $(1-\ell,1)$ for each fixed
$A>0$. Write $\lambda=\lambda(z)$ below. On $I_\lambda$,
$\lambda<1-\rho<2\lambda$ and
$\lambda<1-\rho+z<3\lambda$. By \eqref{eq:ball-assumptions},
\begin{equation}\label{eq:unified-window-bounds}
 w_B\asymp\lambda^{s-1},\qquad
 \frac{|\nabla w_B|}{w_B}\leq C\lambda^{-1},\qquad
 |\varphi|\leq C\lambda^{-1},\qquad
 |\varphi_\rho|+|\varphi_z|\leq C\lambda^{-2}.
\end{equation}
Notice that, by definition, $\int\varphi(\rho,z)\,d\rho=1$.
Since $\lambda'(z)=1$, direct differentiation gives
\begin{equation}\label{eq:unified-window-identity}
 \varphi_z
 =-\frac1{\lambda^2}\chi\!\left(\frac{1-\rho}{\lambda}\right)
  -\frac{1-\rho}{\lambda^3}\chi'\!\left(\frac{1-\rho}{\lambda}\right)
 =\partial_\rho\!\left(\frac{1-\rho}{\lambda}\varphi\right).
\end{equation}

We first explain how to differentiate $H$ away from the bottom.
As noted in Remark~\ref{rem:vertic-deriv}, we may write \eqref{eq:weak-flux} as
$\partial_z(wq_A/w_B)=-\partial_\rho B_A+C_A$, where $B_A,C_A$ are given by \eqref{def:BA-CA}.
 Testing against
$\varphi(\rho,z)\theta(z)$, with $\theta\in C_c^\infty((0,r_0))$, yields in the sense of distributions
\begin{equation}\label{eq:Hphi'}
 H'(z)=\int_\rho\left(B_A\varphi_\rho+C_A\varphi
              +\frac{wq_A}{w_B}\varphi_z\right)d\rho.
\end{equation}
On every compact height interval in $(0,r_0)$ the right-hand side is
integrable: $q_A,(q_A)_\rho$ are bounded, $v\in H^1_{\rm loc}(\mathcal R)$,
and the coefficients are smooth on the horizontal support.
Thus $H\in W^{1,1}_{\rm loc}((0,r_0))$.

Using $w=w_B-v_\rho$ and \eqref{eq:unified-window-identity}, horizontal
integration by parts gives
\[
 \begin{aligned}
 \int\frac{wq_A}{w_B}\varphi_z\,d\rho
 &=-\int\frac{1-\rho}{\lambda}(q_A)_\rho\varphi\,d\rho
   -\int\frac{v_\rho q_A}{w_B}
          \partial_\rho\!\left(\frac{1-\rho}{\lambda}\varphi\right)d\rho\\
 &=-\int\frac{1-\rho}{\lambda}\frac{w}{w_B}(q_A)_\rho\varphi\,d\rho
   -\int\frac{v_\rho}{w_B}
          \partial_\rho\!\left(\frac{1-\rho}{\lambda}q_A\varphi\right)d\rho,
 \end{aligned}
\]
where the second equality uses $1=w/w_B+v_\rho/w_B$.
Substituting in \eqref{eq:Hphi'} gives the exact identity
\begin{equation}\label{eq:Hphi'-expanded2}
 \begin{aligned}
 H'(z)={}&\int_{I_\lambda}
 \left[\frac{J_A}{w_B}
 -\frac{1-\rho}{\lambda}\frac{w}{w_B}(q_A)_\rho\right]\varphi\,d\rho\\
 &+\int_{I_\lambda}v_z\,\partial_\rho\!\left(
                  \frac{q_A\varphi}{w_B}\right)d\rho
 +\int_{I_\lambda}\frac{v_\rho q_A(w_B)_z}{w_B^2}\varphi\,d\rho\\
 &-\int_{I_\lambda}\frac{v_\rho}{w_B}
       \partial_\rho\!\left(\frac{1-\rho}{\lambda}q_A\varphi\right)d\rho.
 \end{aligned}
\end{equation}

We estimate this derivative from below. Since $w\geq0$,
$(1-\rho)/\lambda>0$, $\varphi\geq0$, and
$((q_A)_\rho)_+\leq(q_\rho)_+$, the first line of
\eqref{eq:Hphi'-expanded2} is bounded below by
\begin{equation}\label{Hphi'1}
 -C\lambda^{1/2-s}
 \left[\|J_A(\cdot,z)\|_{L^2(I_\lambda)}
          +\|w(q_\rho)_+(\cdot,z)\|_{L^2(I_\lambda)}\right].
\end{equation}
Indeed, $\varphi/w_B\leq C\lambda^{-s}$,
$(1-\rho)/\lambda\leq2$, and $|I_\lambda|=\lambda$.
For the remaining terms, the product rule, $0\leq q_A\leq A$,
$|(q_A)_\rho|\leq \Lambda$, and \eqref{eq:unified-window-bounds} give
\begin{equation}\label{eq:unified-coefficient-bounds}
 \begin{aligned}
 \left|\partial_\rho\!\left(\frac{q_A\varphi}{w_B}\right)\right|
 &\leq C\left(\Lambda\lambda^{-s}+A\lambda^{-1-s}\right),\\
 \frac{q_A|(w_B)_z|}{w_B^2}|\varphi|
 &\leq CA\lambda^{-1-s},\\
 \frac1{w_B}\left|\partial_\rho\!\left(
           \frac{1-\rho}{\lambda}q_A\varphi\right)\right|
 &\leq C\left(\Lambda\lambda^{-s}+A\lambda^{-1-s}\right).
 \end{aligned}
\end{equation}
For example, in the first line the derivative of $q_A$ contributes
$C\Lambda\lambda^{-s}$, while the derivatives of $\varphi$ and $1/w_B$
each contribute $CA\lambda^{-1-s}$. In the last line,
$|\partial_\rho((1-\rho)/\lambda)|=\lambda^{-1}$ gives the same power.
Cauchy--Schwarz on $I_\lambda$ therefore yields
\[
 H'(z)\geq-C\lambda^{1/2-s}\left[
 \|J_A(\cdot,z)\|_{L^2(I_\lambda)}
 +\|w(q_\rho)_+(\cdot,z)\|_{L^2(I_\lambda)}
 +\left(\Lambda+\frac A\lambda\right)
       \|\nabla v(\cdot,z)\|_{L^2(I_\lambda)}\right].
\]
Since $\lambda\geq\tau=A/K$, we have $A/\lambda\leq K$, and obtain
for almost every $z\in(0,r_0)$
\begin{equation}\label{eq:unified-lower-derivative}
 H'(z)\geq-C\lambda(z)^{1/2-s}\left[
 \|J_A(\cdot,z)\|_{L^2(I_{\lambda(z)})}
 +\|w(q_\rho)_+(\cdot,z)\|_{L^2(I_{\lambda(z)})}
 +\|\nabla v(\cdot,z)\|_{L^2(I_{\lambda(z)})}\right].
\end{equation}

We next claim that
\begin{equation}\label{claim:bottom-average}
 \liminf_{z\downarrow0}H(z)
 \geq A-2\Lambda\tau-C\sqrt{\tau\mathscr J}.
\end{equation}
Since $w=w_B-v_\rho$, horizontal integration by parts gives
\begin{equation}\label{eq:unified-bottom-decomposition}
 H(z)=\int_{I_\lambda}q_A\varphi\,d\rho
       +\int_{I_\lambda}v\,\partial_\rho\!\left(
                          \frac{q_A\varphi}{w_B}\right)d\rho.
\end{equation}
By horizontal Lipschitz continuity,
$q_A(\rho,z)\geq\min\{q(1,z)-2\Lambda\lambda(z),A\}$ on $I_\lambda$.
Since $q(1,z)\to A$, $\lambda(z)\to\tau$, and $\int\varphi=1$,
\[
 \liminf_{z\downarrow0}\int_{I_\lambda}q_A\varphi\,d\rho
 \geq A-2\Lambda\tau.
\]
For $z<\tau$, \eqref{eq:unified-coefficient-bounds} and $A/\tau=K$
bound the coefficient in the second term of
\eqref{eq:unified-bottom-decomposition} by $C\tau^{-s}$.
Write $v(\cdot,z)=v_0+[v(\cdot,z)-v_0]$, where $v_0$ is the bottom
trace. The strong $L^2((1-\ell,1))$ convergence of the slices of
$v$ given by Lemma~\ref{lem:hardy} gives
\[
 \left|\int_{I_\lambda}(v-v_0)\,
       \partial_\rho\!\left(\frac{q_A\varphi}{w_B}\right)d\rho\right|
 \leq C\tau^{-s}\sqrt\lambda\,
       \|v(\cdot,z)-v_0\|_{L^2((1-\ell,1))}\longrightarrow0.
\]
Also $I_{\lambda(z)}\to I_\tau=(1-2\tau,1-\tau)$ in measure, so
integrability of $v_0$ implies
$\int_{I_{\lambda(z)}}|v_0|\,d\rho\to\int_{I_\tau}|v_0|\,d\rho$.
The Hardy term in \eqref{eq:J} gives
\[
 \int_{I_\tau}|v_0|\,d\rho
 \leq\left(\int_{I_\tau}
              \frac{|v_0|^2}{(1-\rho)^{2s}}\,d\rho\right)^{1/2}
      \left(\int_{I_\tau}(1-\rho)^{2s}\,d\rho\right)^{1/2}
 \leq C\tau^{s+1/2}\sqrt{\mathscr J}.
\]
Consequently the limsup of the absolute value of the second term in
\eqref{eq:unified-bottom-decomposition} is at most
\[
 C\tau^{-s}\tau^{s+1/2}\sqrt{\mathscr J}
 =C\sqrt{\tau\mathscr J}.
\]
This proves
\eqref{claim:bottom-average}.

We now propagate this information directly to a fixed height.
Local absolute continuity and \eqref{eq:unified-lower-derivative}
give, for $\varepsilon<z$,
\[
\begin{aligned}
 H(\varepsilon)\leq H(z)+C\int_\varepsilon^z
 \lambda(z')^{1/2-s}\Big[
 &\|J_A(\cdot,z')\|_{L^2(I_{\lambda(z')})}
 +\|w(q_\rho)_+(\cdot,z')\|_{L^2(I_{\lambda(z')})}\\
 &+\|\nabla v(\cdot,z')\|_{L^2(I_{\lambda(z')})}\Big]\,dz'
 \leq H(z)+C\sqrt{\mathscr J},
\end{aligned}
\]
where weighted Cauchy--Schwarz and \eqref{eq:J} apply because
\begin{equation}\label{eq:unified-weight-cancellation}
 \int_0^{r_0}\lambda(z)^{1-2s}z^{2s-1}\,dz
 =\int_0^{r_0}\left(\frac{z+\tau}{z}\right)^{1-2s}\,dz
 \leq Cr_0,
\end{equation}
uniformly in $0<\tau\leq r_0/2$. Combining this with
\eqref{claim:bottom-average} yields,
\begin{equation}\label{eq:unified-propagation}
 A-2\Lambda\tau-C\sqrt{\tau\mathscr J}\leq H(z)+C\sqrt{\mathscr J}.
\end{equation}

Finally, for $z\in[r_0/2,r_0]$, we have $\lambda(z)\asymp r_0$,
$\varphi\leq C$, $w_B\asymp1$, and $z^{1-2s}\asymp1$. Since $0\leq q_A\leq q$ and
$w\geq0$,
\[
\begin{aligned}
 \int_{r_0/2}^{r_0}H(z)\,dz
 &\leq C\int_{r_0/2}^{r_0}\int_{I_{\lambda(z)}}wq\,d\rho\,dz\\
 &\leq C\left(\int_{r_0/2}^{r_0}\int_{I_{\lambda(z)}}
                       z^{1-2s}w^2q^2\,d\rho\,dz\right)^{1/2}
 \leq C\sqrt{\mathscr J}.
\end{aligned}
\]
Integrating \eqref{eq:unified-propagation} in $(r_0/2,r_0)$  gives
\[
 \left(1-\frac{2\Lambda}{K}\right)A=A-2\Lambda\tau
 \leq C\sqrt{\mathscr J}+C\sqrt{\tau\mathscr J}
 \leq C\sqrt{\mathscr J}.
\]
Since $K\geq8\Lambda$, this proves \eqref{eq:trace-result} and hence the
theorem. 
\end{proof}

\section{Reduction to bounded sets}\label{sec:reduction}

In this section we remove the assumption $\Omega\subset B_R$ of
Theorem~\ref{thm:half-confined-SV} by a truncation argument. As in
\cite[Lemma~4.1]{GMR}, which follows \cite[Lemma~5.3]{BDV}, a nonsharp
stability estimate provides the initial smallness of the mass of
$\Omega$ outside a ball, and an iteration on the tail volumes then
reduces it to the order of the deficit. The nonsharp input is
Proposition~\ref{prop:nonsharp-saint-venant}, see also \cite{BCV,Kim}. The local
$L^\infty$ tail estimate for the torsion function used in
\cite{BDV,GMR} is replaced by an $L^1$ tail estimate.

\subsection{The exact cutoff error}

\begin{lemma}\label{lem:cutoff}
Let $u\in H^s(\R^n)$ and let $\chi$ be bounded and
Lipschitz. Then
\begin{equation}\label{eq:cutoff-identity}
 \mathcal E_s(\chi u)=\mathcal E_s(u,\chi^2u)
 +\mathcal R_\chi
\end{equation}
where 
\[
\mathcal R_\chi=\frac{c_{n,s}}2\iint
 \frac{u(x)u(y)(\chi(x)-\chi(y))^2}{|x-y|^{n+2s}}\,dx\,dy.
\]
If $u=u_\Omega$, $0\leq\chi\leq1$, and $\chi u\in H_0^s(D)$ for
an open $D\subset\Omega$, then
\begin{equation}\label{eq:cutoff-torsion}
 T_s(\Omega)-T_s(D)
 \leq\int_\Omega(1-\chi)^2u\,dx+\mathcal R_\chi.
\end{equation}
\end{lemma}

\begin{proof}
For real numbers $a,b,c,d$,
$(ca-db)^2-(a-b)(c^2a-d^2b)=ab(c-d)^2$.
Apply this pointwise to the two values of $u$ and $\chi$ and integrate. This gives \eqref{eq:cutoff-identity}. Since $\chi$ is Lipschitz, the error term $\mathcal R_\chi$
 is absolutely integrable.

If $u=u_\Omega$, testing \eqref{eq:torsion-weak}  with $\chi^2u$ gives
$\mathcal E_s(u,\chi^2u)=\int\chi^2u$. Insert $\chi u$ into the
variational formula for $T_s(D)$ to obtain \eqref{eq:cutoff-torsion}.
\end{proof}

\subsection{Truncation}

\begin{proposition}\label{prop:half-truncation}
There exist $R_*,C,\delta_*>0$, depending
only on $n$ and $s$, such that every open set
$\Omega\subset\R^n$ with $|\Omega|=\omega_n$ and
$\delta_T(\Omega)\leq\delta_*$ admits an open set
$\widetilde\Omega\subset B_{R_*}$ satisfying
\begin{equation}\label{eq:truncation-target}
 \begin{gathered}
 |\widetilde\Omega|=\omega_n,\qquad
 \delta_T(\widetilde\Omega)\leq C\delta_T(\Omega),\\
 \mathcal A(\Omega)\leq\mathcal A(\widetilde\Omega)
                         +C\delta_T(\Omega).
 \end{gathered}
\end{equation}
\end{proposition}

\begin{proof}
Write $\delta=\delta_T(\Omega)$ and $u=u_\Omega$. Since $n\geq2$
and $s<1$, we have $n>2s$.

First,
\begin{equation}\label{eq:truncation-global-Linfty}
 \|u\|_{L^\infty(\R^n)}\leq C.
\end{equation}
To prove this, set $2_s^*=2n/(n-2s)$ and
$\mu(t)=|\{u>t\}|$. Let $f_t=(u-t)_+$. Since for all real $a,b$,
\[(a-b)\bigl((a-t)_+-(b-t)_+\bigr)\geq\bigl((a-t)_+-(b-t)_+\bigr)^2,\]
 testing \eqref{eq:torsion-weak} with $f_t$ gives
$\mathcal E_s(f_t)\leq\mathcal E_s(u,f_t)=\int_\Omega f_t$. Combining this
with fractional Sobolev and H\"older's inequality (recall that $f_t$ is
supported in $\{u>t\}$), we get
\[
 \|f_t\|_{2_s^*}^2\leq C\mathcal E_s(f_t)
 \leq C\int_\Omega f_t
 \leq C\mu(t)^{1-1/2_s^*}\|f_t\|_{2_s^*}.
\]
Consequently, for $t'>t\geq0$,
\[
 \mu(t')\leq
 \left(\frac{C}{t'-t}\right)^{2_s^*}\mu(t)^{2_s^*-1}.
\]
With $t_j=\Theta(1-2^{-j})$ and
$\Theta=K\mu(0)^{2s/n}$, induction yields
$\mu(t_j)\leq2^{-nj/(2s)}\mu(0)$ if $K$ is sufficiently large.

Suppose $0\leq\chi\leq1$, $|\nabla\chi|\leq C_1$, and $\chi$ is
constant on $B_k$, where $C_1$ is fixed independently of $k$.
If $\chi(x)\ne\chi(y)$, at least one of $x,y$ lies outside $B_k$.
Since $u\geq0$, symmetry and \eqref{eq:truncation-global-Linfty} give
\begin{equation}\label{eq:truncation-error-tail}
 \begin{aligned}
 0\leq\mathcal R_\chi
 &\leq C\int_{\Omega\setminus B_k}u(x)
       \int_{\R^n}u(y)
       \frac{\min\{1,C_1^2|x-y|^2\}}{|x-y|^{n+2s}}\,dy\,dx\\
 &\leq C\int_{\Omega\setminus B_k}u.
 \end{aligned}
\end{equation}

For $k\geq1$, put $S_k=\Omega\setminus B_k$ and choose a smooth
cutoff $\chi_k$ such that
\[
 0\leq\chi_k\leq1,\qquad
 \chi_k=0\text{ on }B_k,\qquad
 \chi_k=1\text{ outside }B_{k+1},\qquad
 |\nabla\chi_k|\leq C.
\]
 By
\eqref{eq:cutoff-identity} and \eqref{eq:truncation-error-tail},
\begin{equation}\label{eq:truncation-exterior-energy}
 \mathcal E_s(\chi_ku)
 =\int_\Omega\chi_k^2u+\mathcal R_{\chi_k}
 \leq C\int_{S_k}u\leq C|S_k|.
\end{equation}
Since $\chi_ku=u$ on $S_{k+1}$, Sobolev and H\"older imply
\begin{equation}\label{eq:truncation-tail-L1}
 \begin{aligned}
 \int_{S_{k+1}}u
 &\leq |S_{k+1}|^{1-1/2_s^*}\|\chi_ku\|_{2_s^*}\\
 &\leq C|S_k|^{1/2+s/n}
                     \mathcal E_s(\chi_ku)^{1/2}
 \leq C|S_k|^{1+s/n}.
 \end{aligned}
\end{equation}

Next choose a smooth cutoff $\zeta_k$ satisfying
\[
 0\leq\zeta_k\leq1,\qquad
 \zeta_k=1\text{ on }B_{k+1},\qquad
 \operatorname{supp}\zeta_k\Subset B_{k+2},\qquad
 |\nabla\zeta_k|\leq C.
\]
For $D_k=\Omega\cap B_{k+2}$, applying \eqref{eq:cutoff-torsion},
\eqref{eq:truncation-error-tail} with inner radius $k+1$, and
\eqref{eq:truncation-tail-L1}, we obtain
\begin{equation}\label{eq:truncation-torsion-loss}
 \begin{aligned}
 0\leq T_s(\Omega)-T_s(D_k)
 &\leq\int_\Omega(1-\zeta_k)^2u+\mathcal R_{\zeta_k}\\
 &\leq C\int_{S_{k+1}}u
 \leq C|S_k|^{1+s/n}.
 \end{aligned}
\end{equation}

Proposition~\ref{prop:nonsharp-saint-venant} gives
$\mathcal A(\Omega)\leq C\delta^{s/(2s+1)}$.
If $\delta=0$, then $\mathcal A(\Omega)=0$ and we take
$\widetilde\Omega=B$. Otherwise translate $\Omega$ so that $B_1$
realizes its Fraenkel asymmetry.  Set
\[
 a_k=\frac{|\Omega\setminus B_k|}{\omega_n},\qquad
 \sigma=1+\frac{s}{n}.
\]
Then
\begin{equation}\label{eq:truncation-seed}
 a_1=\tfrac12\mathcal A(\Omega)\leq C\delta^{s/(2s+1)}.
\end{equation}
The fractional Saint-Venant inequality $T_s(D_k)\leq T_s(B_{r_k})$,
where $|B_{r_k}|=|D_k|$, which follows from Lemma~\ref{lem:split} and
\eqref{eq:scaling}, together with
\eqref{eq:truncation-torsion-loss}, give
\[
 T_s(B)-\delta-Ca_k^\sigma
 \leq T_s(D_k)
 \leq T_s(B)(1-a_{k+2})^{1+2s/n}.
\]
Since $1-(1-a)^{1+2s/n}\geq a$ for $a\in[0,1]$, it follows that
\begin{equation}\label{eq:truncation-recurrence}
 a_{k+2}\leq C\bigl(\delta+a_k^\sigma\bigr).
\end{equation}
For $0<\delta\leq1$, induction gives
\begin{equation}\label{eq:truncation-iteration}
 a_{2j+1}\leq C_j\delta^{\min\{s\sigma^j/(2s+1),1\}}.
\end{equation}
 Choose $N$ with
$s\sigma^N/(2s+1)\geq1$ and
put $k_*=2N+1$. Then
\begin{equation}\label{eq:truncation-small-tail}
 a_{k_*}\leq C_N\delta.
\end{equation}

Let
\[
 D=\Omega\cap B_{k_*+2},\qquad
 a=\frac{|\Omega\setminus D|}{\omega_n},\qquad
 r=(1-a)^{1/n}.
\]
Choose $\delta_*\leq\min\{1,(2C_N)^{-1}\}$. Then
$a\leq a_{k_*}\leq1/2$, and $\widetilde\Omega=r^{-1}D$ satisfies
\[
 |\widetilde\Omega|=\omega_n,\qquad
 \widetilde\Omega\subset B_{R_*},\qquad
 R_*=2^{1/n}(k_*+2).
\]
By \eqref{eq:truncation-torsion-loss} and
\eqref{eq:truncation-small-tail},
\[
 T_s(\Omega)-T_s(D)
 \leq Ca_{k_*}^\sigma\leq C\delta.
\]
Since $T_s(\widetilde\Omega)=r^{-(n+2s)}T_s(D)
\geq T_s(D)$, this yields
\[
 \delta_T(\widetilde\Omega)
 \leq\delta+
 T_s(\Omega)-T_s(D)
 \leq C\delta.
\]
If $B_1(b)$ realizes the asymmetry of $\widetilde\Omega$, comparison
with $B_1(rb)$ gives
\[
 \begin{aligned}
 \omega_n\mathcal A(\Omega)
 &\leq|\Omega\setminus D|+|D\triangle B_r(rb)|
                         +|B_r(rb)\triangle B_1(rb)|\\
 &=2\omega_n a+r^n\omega_n\mathcal A(\widetilde\Omega).
 \end{aligned}
\]
The last assertion of \eqref{eq:truncation-target} follows from
$a\leq C_N\delta$ and $r\leq1$.
\end{proof}

\begin{proof}[Proof of Theorem~\ref{thm:main-sv}]
If $\delta_T(\Omega)\leq\delta_*$, apply
Proposition~\ref{prop:half-truncation} and then
Theorem~\ref{thm:half-confined-SV} to $\widetilde\Omega$ to obtain
\[
 \mathcal A(\Omega)^2
 \leq C\bigl(\mathcal A(\widetilde\Omega)^2
                    +\delta_T(\Omega)^2\bigr)
 \leq C\delta_T(\Omega).
\]
If $\delta_T(\Omega)>\delta_*$, the conclusion follows from
$\mathcal A(\Omega)\leq2$.
\end{proof}

\section{Optimality of the exponent}\label{sec:eigenvalue}

We now show that the exponent 2  in Theorems~\ref{thm:main-fk} and~\ref{thm:main-sv} is indeed optimal. As in \cite{BDV} it can be shown by estimating the deficits on ellipsoids.

\begin{lemma}\label{lem:ellipsoidal-sharpness}
Let
\[
 L_\eps=\operatorname{diag}(e^\eps,e^{-\eps},1,\ldots,1),
 \qquad \Omega_\eps=L_\eps B.
\]
For every fixed $n\geq2$ and $s\in(0,1)$, these sets satisfy
$|\Omega_\eps|=\omega_n$ and, for all sufficiently small $|\eps|$,
\begin{equation}\label{eq:ellipsoidal-sharpness}
 \begin{gathered}
 c|\eps|\leq\mathcal A(\Omega_\eps)\leq C|\eps|,\\
 0\leq\frac{\lambda_s(\Omega_\eps)}{\lambda_s(B)}-1
 \leq M_s(\eps)-1,\qquad
 0\leq\frac{T_s(B)-T_s(\Omega_\eps)}{T_s(B)}
 \leq1-\frac1{M_s(\eps)},
 \end{gathered}
\end{equation}
where
\begin{equation}\label{eq:ellipsoidal-angular-factor}
 \begin{aligned}
 M_s(\eps)&=\frac1{|\mathbb S^{n-1}|}
 \int_{\mathbb S^{n-1}}|L_\eps^{-T}\theta|^{2s}\,d\mathcal H^{n-1}(\theta)\\
 &=1+\frac{4s(n+2s)}{n(n+2)}\eps^2+O(\eps^4).
 \end{aligned}
\end{equation}
The constants and the remainder bound depend only on $n,s$.
\end{lemma}

\begin{proof}
The volume assertion follows from $\det L_\eps=1$.
For any radial $f\in H_0^s(B)$, put
$f_\eps(x)=f(L_\eps^{-1}x)$. Linear changes of variables in smooth
compactly supported approximations show that
$f_\eps\in H_0^s(\Omega_\eps)$. Its integral and $L^2$ norm agree
with those of $f$. Moreover,
$\widehat{f_\eps}(\xi)=\widehat f(L_\eps^T\xi)$, so the Fourier
representation of the quadratic form and polar integration give
\begin{equation}\label{eq:ellipsoidal-energy}
 \begin{aligned}
 \mathcal E_s(f_\eps)
 &=\frac1{(2\pi)^n}\int_{\R^n}
 |L_\eps^{-T}\xi|^{2s}|\widehat f(\xi)|^2\,d\xi\\
 &=M_s(\eps)\mathcal E_s(f).
 \end{aligned}
\end{equation}
Here $\widehat f$ is radial. A nonnegative radial first eigenfunction
$\phi$ of $B$ can be chosen by Schwarz rearrangement of a minimizer
of the Rayleigh quotient. Applying \eqref{eq:ellipsoidal-energy} to
$\phi$ gives
$\lambda_s(\Omega_\eps)\leq M_s(\eps)\lambda_s(B)$.
For torsion, take $f=u_B$ and use the trial functions
$a f_\eps$, $a\in\R$, in \eqref{eq:torsion-rigidity}. Since
$\int f_\eps=\mathcal E_s(u_B)=T_s(B)$,
\[
 T_s(\Omega_\eps)\geq
 \max_{a\in\R}\bigl(2a-a^2M_s(\eps)\bigr)T_s(B)
 =\frac{T_s(B)}{M_s(\eps)}.
\]
The qualitative Faber--Krahn and Saint-Venant inequalities give the
nonnegativity of both deficits, and also $M_s(\eps)\geq1$.

To compute the angular factor, write
$a(\theta)=\theta_1^2-\theta_2^2$ and
$b(\theta)=\theta_1^2+\theta_2^2$. Uniformly on the sphere,
\[
 |L_\eps^{-T}\theta|^{2s}
 =1-2s\eps a(\theta)
 +\eps^2\bigl(2s b(\theta)+2s(s-1)a(\theta)^2\bigr)
 +O(|\eps|^3).
\]
The spherical averages of $a$, $b$, and $a^2$ are respectively
$0$, $2/n$, and $4/[n(n+2)]$. This gives the quadratic coefficient
in \eqref{eq:ellipsoidal-angular-factor}. Interchanging the first two
coordinates shows that $M_s$ is even; its smoothness near zero and
Taylor's theorem improve the remainder to $O(\eps^4)$.

It remains to estimate the asymmetry, including the infimum over
translations. The radial function of $\Omega_\eps$ is
\[
 r_\eps(\theta)=|L_\eps^{-1}\theta|^{-1}
 =1+\eps(\theta_1^2-\theta_2^2)+O(\eps^2),
\]
with uniform remainder. Consequently,
\[
 \begin{aligned}
 |\Omega_\eps\triangle B|
 &=\frac1n\int_{\mathbb S^{n-1}}|r_\eps(\theta)^n-1|
 \,d\mathcal H^{n-1}(\theta)\\
 &=|\eps|\int_{\mathbb S^{n-1}}|\theta_1^2-\theta_2^2|
 \,d\mathcal H^{n-1}(\theta)+O(\eps^2).
 \end{aligned}
\]
The coefficient is positive because $n\geq2$.
Since $\Omega_\eps=-\Omega_\eps$, for every center $a\in\R^n$,
\[
 |B_1(a)\triangle B_1(-a)|
 \leq2|\Omega_\eps\triangle B_1(a)|.
\]
The symmetric difference of two unit balls is nondecreasing in the
distance between their centers, as follows by slicing parallel to
that displacement. Thus
\[
 \begin{aligned}
 |\Omega_\eps\triangle B|
 &\leq|\Omega_\eps\triangle B_1(a)|+|B_1(a)\triangle B|\\
 &\leq3|\Omega_\eps\triangle B_1(a)|.
 \end{aligned}
\]
Taking the infimum over $a$ proves
\[
 \frac{|\Omega_\eps\triangle B|}{3\omega_n}
 \leq\mathcal A(\Omega_\eps)
 \leq\frac{|\Omega_\eps\triangle B|}{\omega_n},
\]
which completes the proof.
\end{proof}

\appendix

\section{Ball speed and boundary control of the radial profile}
\label{app:ball-hardy}

\subsection{Estimates for the ball speed}

\begin{lemma}\label{lem:half-ball-speed}
Let $U_B=\operatorname{Ext}_s u_B$, with
$w_B=-\partial_\rho t_B$. There are $\ell\in(0,1/16)$ and constants
$c_0,C_0>0$, depending only on $n,s$, such that the ball speed is
positive on $\mathcal R$ and
\begin{equation}\label{eq:ball-assumptions}
 \begin{gathered}
 c_0(1-\rho+z)^{s-1}\leq w_B(\rho,z)
 \leq C_0(1-\rho+z)^{s-1},\\
 |\nabla w_B(\rho,z)|\leq C_0\frac{w_B(\rho,z)}{1-\rho+z}.
 \end{gathered}
\end{equation}
Moreover, $w_B$ is smooth up to every compact subinterval of the bottom
edge and
$(w_B)_z(\rho,0)=0$ for $0<\rho<1$.
\end{lemma}

\begin{proof}
The explicit ball formula \cite[Theorem~1]{Dyda}, with exponent $s$
and fractional order $2s$, gives in our normalization
\begin{equation}\label{eq:general-ball-torsion}
 u_B(x)=\kappa_{n,s}(1-|x|^2)_+^s,\qquad
 \kappa_{n,s}=\frac{\Gamma(n/2)}
 {2^{2s}\Gamma(1+s)\Gamma(n/2+s)}.
\end{equation}
The operator in \cite{Dyda} has the opposite sign; reversing it gives
the positive multiplier $|\xi|^{2s}$. In particular,
\begin{equation}\label{eq:general-ball-bottom-speed}
 w_B(\rho,0)=2s\kappa_{n,s}\rho(1-\rho^2)^{s-1}
 \asymp(1-\rho)^{s-1}\qquad(1/2<\rho<1).
\end{equation}
We justify below that this radial derivative of the trace is also
the bottom trace of $w_B$.

\emph{Weighted local estimates and even reflection.}
We use the local boundedness, H\"older continuity, and Harnack estimates
for weak solutions of
\[
 \operatorname{div}(|z|^{1-2s}\nabla f)=0
\]
in \cite[Section~2.3]{FKS}. The weight $|z|^{1-2s}$ belongs to $A_2$
precisely because $-1<1-2s<1$. The constants depend on the dimension
and $s$, and the estimates are invariant under simultaneous dilation
of the spatial variables and height. Harnack applies to nonnegative
solutions, including balls that meet $z=0$.

We also need tangential derivative estimates, which follow from these
inputs as follows. The coefficient is independent of the horizontal
variables, so horizontal difference quotients solve the same equation.
The weighted Caccioppoli inequality, iterated on nested balls, bounds
the weighted $L^2$ norms of all horizontal derivatives. Local
boundedness applied to those derivatives then gives, on a smaller ball
of radius comparable to $r$,
\begin{equation}\label{eq:weighted-tangential-estimates}
 \sup |D_y^m f|\leq C_m r^{-m}\sup |f|\qquad(m\geq1),
\end{equation}
where the supremum on the right is over a fixed larger ball and $y$
denotes the horizontal variables in the dimension being considered.
For signed derivatives one applies local boundedness to their positive
and negative parts. H\"older continuity gives continuous representatives
of all these horizontal derivatives across $z=0$.

If $f$ is even in $z$ and is a weak solution across $z=0$, its weighted
normal flux has zero trace. Indeed,
$\partial_z(|z|^{1-2s}f_z)=-|z|^{1-2s}\Delta_Yf$ has an integrable
right-hand side locally, by the preceding tangential regularity.
The flux is continuous across zero and odd, so its trace there is zero.
Integration in $z>0$ therefore gives the exact identity
\begin{equation}\label{eq:weighted-even-normal-derivative}
 f_z(y,z)=-z^{2s-1}\int_0^z t^{1-2s}\Delta_Yf(y,t)\,dt
        =-z\int_0^1 t^{1-2s}\Delta_Yf(y,tz)\,dt.
\end{equation}
In particular $f_z(y,0)=0$, and
\begin{equation}\label{eq:weighted-even-normal-bound}
 |f_z(y,z)|\leq\frac{z}{2-2s}
                  \sup_{0\leq t\leq z}|\Delta_Yf(y,t)|.
\end{equation}
The last expression in \eqref{eq:weighted-even-normal-derivative},
first for all horizontal derivatives and then differentiated repeatedly
in $z$, proves smoothness across the bottom. More precisely, continuity
of all horizontal derivatives starts the induction on the number of
normal derivatives; at each step differentiation under the integral
uses an integrable power $t^{1-2s+j}$, $j\geq0$. The even extension
satisfies the same formula for negative $z$.

\emph{The ball equation and its zero-flux derivative.}
On each bottom patch compactly contained in $B$, the weak torsion
condition is
\[
 -d_s\lim_{z\downarrow0}z^{1-2s}\partial_zU_B=1.
\]
Since $\partial_z(z^{2s}/(2s d_s))=z^{2s-1}/d_s$, the function
$U_B+z^{2s}/(2s d_s)$ has zero weak weighted normal flux on this
patch. Its even reflection is a weak solution of
$\operatorname{div}(|z|^{1-2s}\nabla f)=0$ across the patch, by
\cite[Lemma~4.1]{CS}, or directly by testing on the two half-balls.
The added term has finite local weighted energy because
$\int_0^r z^{2s-1}\,dz<\infty$. The preceding local argument shows
that this corrected reflection is smooth. Its radial derivative is
$-w_B$, since the correction is independent of the horizontal variables.
Thus $w_B$ has a smooth even reflection, its bottom value is
\eqref{eq:general-ball-bottom-speed}, and $(w_B)_z(\rho,0)=0$.

The convolution of two nonnegative radial nonincreasing functions is
radial nonincreasing, so $w_B\geq0$ at positive height. Differentiating
the weighted radial equation for $t_B$ gives
\[
 (w_B)_{\rho\rho}+\frac{n-1}\rho(w_B)_\rho+(w_B)_{zz}
 +\frac{1-2s}{z}(w_B)_z-\frac{n-1}{\rho^2}w_B=0.
\]
Equivalently, $k(\rho,z)=w_B(\rho,z)/\rho$ satisfies
\begin{equation}\label{eq:general-ball-lifted-equation}
 k_{\rho\rho}+\frac{n+1}\rho k_\rho+k_{zz}
 +\frac{1-2s}{z}k_z=0.
\end{equation}
Thus $k(|y|,z)$, for $y\in\R^{n+2}$, solves
$\operatorname{div}(z^{1-2s}\nabla k)=0$ away from $|y|=0$.
It is nonnegative and not identically zero, hence positive at positive
height by the strong maximum principle. On bottom annuli with
$0<|y|<1$, its smooth even reflection solves the corresponding equation
with weight $|z|^{1-2s}$: its normal derivative is $O(z)$, so the
weighted flux is $O(z^{2-2s})\to0$ and there is no interface source.
All subsequent estimates stay away from $|y|=0$.

\emph{The regime $0<z\leq1-\rho$.}
Write $d=1-\rho$ and $y_0=\rho e_1\in\R^{n+2}$. For $d$ small,
the cylinder $|y-y_0|<d/2$, $|z|<2d$ lies in the reflected region
and has $|y|>1/2$. A bounded Harnack chain of balls with radius
comparable to $d$ joins $(y_0,0)$ to $(y_0,z)$ and covers a smaller
neighborhood of the segment between them. The weighted Harnack
inequality and \eqref{eq:general-ball-bottom-speed} imply
\begin{equation}\label{eq:general-ball-low-height}
 k(\rho,z)\asymp k(\rho,0)\asymp d^{s-1}
 \qquad(0<z\leq d).
\end{equation}
The number of balls and their relative radii are independent of $d,z$.
Applying \eqref{eq:weighted-tangential-estimates} for $m=1,2$ on
these balls and using Harnack to compare their suprema with $k(\rho,z)$
gives
\[
 |k_\rho(\rho,z)|\leq C d^{-1}k(\rho,z),\qquad
 \sup_{0\leq t\leq z}|\Delta_Yk(y_0,t)|\leq C d^{-2}k(\rho,z).
\]
Consequently \eqref{eq:weighted-even-normal-bound} gives
\[
 |k_z(\rho,z)|\leq C z d^{-2}k(\rho,z)
                     \leq C d^{-1}k(\rho,z).
\]
Since $w_B=\rho k$, $\rho\asymp1$, and $d+z\asymp d$, this proves
\eqref{eq:ball-assumptions} in this regime, after choosing $\ell$ small.

\emph{The regime $z>1-\rho$.}
Join $(\rho e_1,z)$ to $((1-z)e_1,z)$ by a bounded Harnack chain
of balls of radius comparable to $z$. These balls remain in the upper
half-space and have radial coordinates greater than $1/2$. At the
second endpoint \eqref{eq:general-ball-low-height} applies, with bottom
distance $z$, so $k(\rho,z)\asymp z^{s-1}$. On an interior ball of
radius $z/8$ about the first endpoint, scaling by $z$ makes the
coefficient in \eqref{eq:general-ball-lifted-equation} smooth and
uniformly elliptic, with bounds depending only on $n,s$.
The interior gradient estimate and Harnack give
$|\nabla k(\rho,z)|\leq Ck(\rho,z)/z$ there. Multiplying by
$\rho\asymp1$ and using $1-\rho+z\asymp z$ proves the remaining
estimates in \eqref{eq:ball-assumptions}.
\end{proof}

\subsection{Boundary control of the radial profile}

\begin{lemma}\label{lem:hardy}
For every $f\in H^s(\R^n)$ vanishing a.e. outside $B$,
\begin{equation}\label{eq:exterior-hardy}
 \int_B\frac{|f(x)|^2}{(1-|x|)^{2s}}\,dx
 \leq C\,\mathcal E_s(f).
\end{equation}
Let now $\Omega\subset\R^n$ be bounded and open, with
$|\Omega|=\omega_n$, let $u=u_\Omega$, let $0<\ell<1/4$, and consider the
profile error
$v=t-t_B$ of \eqref{eq:corner-notation}, with $v(\cdot,0)$ defined as
the radial profile of $u^*-u_B$. Then
\begin{equation}\label{eq:profile-budget}
 \int_{\mathcal R}z^{1-2s}|\nabla v|^2\,d\rho\,dz
 +\int_{1-\ell}^1\frac{|v(\rho,0)|^2}{(1-\rho)^{2s}}\,d\rho
 \leq C(\ell)\delta_T(\Omega).
\end{equation}
Moreover $v(\cdot,z)\to v(\cdot,0)$ in $L^2((1-\ell,1))$ as
$z\downarrow0$.
\end{lemma}

\begin{proof}
\emph{The exterior Hardy inequality.} Since $f$ vanishes a.e. outside
$B$, the integrand of $\mathcal E_s(f)$ equals
$|f(x)|^2|x-y|^{-n-2s}$ on $B\times B^c$ and
$|f(y)|^2|x-y|^{-n-2s}$ on $B^c\times B$. Discarding the remaining
nonnegative contribution of $B\times B$ and using the symmetry of these
two regions,
\[
 \mathcal E_s(f)\geq\frac{c_{n,s}}2\cdot2\int_B |f(x)|^2
       \int_{B^c}|x-y|^{-n-2s}\,dy\,dx
 =c_{n,s}\int_B |f(x)|^2
       \int_{B^c}|x-y|^{-n-2s}\,dy\,dx.
\]
Put $d=1-|x|\in(0,1]$. If $x\ne0$, let $e=x/|x|$; if $x=0$, take any
unit vector. Every point of the ball $B_{d/2}((1+d)e)$ is at distance at
least $1+d-d/2>1$ from the origin, so that ball lies outside $B$; since
$|x-(1+d)e|=2d$, its points are at distance at most $5d/2$ from $x$.
Hence the inner integral is at least
\[
 \omega_n\left(\frac d2\right)^n\left(\frac{5d}2\right)^{-n-2s}
 =\omega_n2^{2s}5^{-n-2s}d^{-2s},
\]
which proves \eqref{eq:exterior-hardy} with
$C=\bigl(c_{n,s}\omega_n2^{2s}5^{-n-2s}\bigr)^{-1}$.

\emph{The boundary term.} Apply \eqref{eq:exterior-hardy} to
$f=u^*-u_B$. This is admissible: $u_B\in H_0^s(B)$, and
Lemma~\ref{lem:split} shows $u^*\in H_0^s(B)$; in particular both belong
to $H^s(\R^n)$ and vanish a.e. outside $B$, the latter because
$|\{u^*>0\}|\leq|\Omega|=|B|$ and $u^*$ is radial nonincreasing.
Both functions being radial, $f(\rho e_1)=v(\rho,0)$, and polar
integration gives
\[
 \int_0^1P_\rho\frac{|v(\rho,0)|^2}{(1-\rho)^{2s}}\,d\rho
 =\int_B\frac{|f(x)|^2}{(1-|x|)^{2s}}\,dx
 \leq C\mathcal E_s(u^*-u_B).
\]
By \eqref{eq:torsion-split} and the nonnegativity of
$\mathcal E_s(u)-\mathcal E_s(u^*)$ we have
$\mathcal E_s(u^*-u_B)\leq\delta_T(\Omega)$. On $\rho\in[1-\ell,1]$, with
$\ell<1/4$, the radial factor $P_\rho=n\omega_n\rho^{n-1}$ is bounded
above and below by positive dimensional constants. This proves the
boundary part of \eqref{eq:profile-budget}.

\emph{The interior term.} Both $U^*$ and $U_B$ are radial in $x$, so
$U^*-U_B$ has radial profile $v$ and, as in the proof of
Lemma~\ref{lem:half-four-deficits},
\[
 \int_{\R^{n+1}_+}z^{1-2s}|\nabla(U^*-U_B)|^2\,dx\,dz
 =\int_0^\infty z^{1-2s}\int_0^\infty P_\rho|\nabla v|^2\,d\rho\,dz .
\]
Restricting the nonnegative integrand to $\mathcal R$, bounding $P_\rho$
from below there, and using \eqref{eq:split-weighted-profile-control}
gives
\[
 \int_{\mathcal R}z^{1-2s}|\nabla v|^2\,d\rho\,dz
 \leq C\delta_T(\Omega) .
\]

\emph{Identification of the bottom trace.} The kernel $P_z^{(s)}$ is an
approximate identity, so $U_B(\cdot,z)\to u_B$ and $U(\cdot,z)\to u$ in
$L^2(\R^n)$, and \eqref{eq:rearranged-trace} gives
$\|U^*(\cdot,z)-u^*\|_2\leq\|U(\cdot,z)-u\|_2\to0$. Hence
$\int_0^\infty P_\rho|v(\rho,z)-v(\rho,0)|^2\,d\rho\to0$, and the two-sided
bound for $P_\rho$ gives the asserted $L^2((1-\ell,1))$
convergence.
\end{proof}

\section*{Acknowledgments}
The authors acknowledge financial support from the ANR grant STOIQUES.
MN and BR are members of INdAM/GNAMPA. BR is a member of INdAM--GNAMPA and has been partially supported by the ``INdAM - GNAMPA Project'', code CUP E53C25002010001.

\bibliographystyle{acm}

	\bibliography{biblio}

\end{document}